\documentclass{amsart}

\usepackage{graphicx,psfrag,amssymb}

\newcommand{\HH}[3]{\ensuremath{\mathcal{H}_{#1}^{#2,#3}}}
\newcommand{\Hpqr}{\HH{p}{q}{r}}
\newcommand{\HEx}{\HH{3}{5}{2}}
\newcommand{\DD}[2]{\ensuremath{\mathcal{D}_{#1,#2}}}
\newcommand{\Dqr}{\DD{q}{r}}
\newcommand{\A}[2]{\ensuremath{\mathcal{A}_{#1}^{#2}}}
\newcommand{\Apr}{\A{p}{r}}
\newcommand{\C}[1][\alpha]{\ensuremath{\mathcal{C}(#1)}}
\newcommand{\eps}{\varepsilon}
\newcommand{\pent}{\ensuremath{\mathcal{P}}}

\newcommand{\integer}{\ensuremath{\mathbb{Z}}}
\newcommand{\Z}[1]{\ensuremath{\integer/#1\integer}}
\newcommand{\F}{\mathbb{F}}
\newcommand{\Fp}{{\mathbb F}_p}
\newcommand{\Fq}{{\mathbb F}_q}
\newcommand{\Fpbar}{\overline{\F}_p}
\newcommand{\PP}{\mathbb{P}}
\newcommand{\pone}[1][p]{\ensuremath{\PP^1(\F_{#1})}}
\newcommand{\D}[1]{\ensuremath{D_{#1}}}

\newcommand{\psl}[1][p]{\ensuremath{PSL(2,#1)}}
\newcommand{\pgl}[1]{\ensuremath{PGL(2,#1)}}
\newcommand{\SL}[1]{\ensuremath{SL(2,#1)}}
\newcommand{\gl}[1]{\ensuremath{GL(2,#1)}}

\DeclareMathOperator{\Stab}{Stab}
\DeclareMathOperator{\ord}{ord}
\DeclareMathOperator{\Hom}{Hom}
\newcommand{\id}{1}

\newcommand{\la}{\langle}
\newcommand{\ra}{\rangle}

\newcommand{\twobytwo}[4]{\begin{bmatrix} #1 &#2 \\ #3 &#4 \end{bmatrix}}

\numberwithin{equation}{section}

\theoremstyle{plain}
\newtheorem{theorem}{Theorem}[section]
\newtheorem{lemma}[theorem]{Lemma}
\newtheorem{proposition}[theorem]{Proposition}

\theoremstyle{remark}
\newtheorem{remark}[theorem]{Remark}

\begin{document}

\title[Generalised knot groups distinguish the square and granny knots]
      {Generalised knot groups distinguish the square and granny knots\\
       (with an appendix by David Savitt)}
\author{Christopher Tuffley}
\address{Institute of Fundamental Sciences, Massey University,
         Private Bag 11 222, Palmerston North 4442, New Zealand}
\email{C.Tuffley@massey.ac.nz}

\begin{abstract}
Given a knot $K$ we may construct a group $G_n(K)$ from the
fundamental group of $K$ by adjoining an $n$th root of the meridian
that commutes with the corresponding longitude.  These ``generalised
knot groups'' were introduced independently by Wada and Kelly, and
contain the fundamental group as a subgroup.  The
square knot $SK$ and the granny knot $GK$ are a well known example of
a pair of distinct knots with isomorphic fundamental groups. We show
that $G_n(SK)$ and $G_n(GK)$ are non-isomorphic for all $n\geq 2$.
This confirms a conjecture of Lin and Nelson, and shows that the
isomorphism type of $G_n(K)$, $n\geq 2$, carries more information
about $K$ than the isomorphism type of the fundamental group.
The appendix contains some results on representations of the trefoil
group in \psl\ that are needed for the
proof.
\end{abstract}



\maketitle

\section{Introduction}

Wada~\cite{wada92} and Kelly~\cite{kelly90} independently introduced a
family of link invariants $G_n(L)$ generalising the fundamental group
of a link $L$ in $S^3$. These groups may be defined via Wirtinger-type
presentations, with conjugation by the generator $x_j$ corresponding
to the over-arc replaced by conjugation by $x_j^n$, but also admit a
description in terms of the peripheral system of $L$. The group
$G_1(L)$ is simply the fundamental group of $L$, and the second
description of $G_n$ shows that $G_1(L)$ is a subgroup of $G_n(L)$ for
each $n$.

The square knot $SK$ and the granny knot $GK$ are a well known example
of a pair of distinct knots with isomorphic fundamental groups. They
are therefore a natural choice for a pair of knots on which to test
the strength of the invariant $G_n$. Lin and
Nelson~\cite{lin-nelson08} report the results of computer experiments
in which they were unable to distinguish $G_n(SK)$ and $G_n(GK)$
by counting homomorphisms to selected finite groups of orders as
large as $360$. Although they were not successful in distinguishing
$G_n(SK)$ and $G_n(GK)$, they nevertheless conjectured that these
groups were non-isomorphic for each $n\geq 2$. The purpose of this
paper is to confirm their conjecture, and thereby show that the
isomorphism types of the generalised knot groups carry more information
about $K$ than the isomorphism type of the fundamental group itself.

We will distinguish $G_n(SK)$ and $G_n(GK)$ for $n\geq 2$ by
comparing the number of homomorphisms into suitably chosen 
finite groups. 
Our target groups will be wreath products over \psl, and are
described in Section~\ref{targets.sec}. 
Our main result is the following theorem:

\begin{theorem}
\label{maintheorem}
For each $n\geq 2$ there is a finite group $H$ such that
\[
|\Hom(G_n(GK),H)| < |\Hom(G_n(SK),H)|.
\]
Consequently, $G_n(GK)$ is not isomorphic to $G_n(SK)$ for $n\geq 2$.
\end{theorem}

\begin{remark}
Although our work confirms
Lin and Nelson's conjecture that
$G_n(SK)$ and $G_n(GK)$ are non-isomorphic for all $n\geq 2$, it also
provides counterexamples to their initial claim that 
$|\Hom(G_n(SK),H)| = |\Hom(G_n(GK),H)|$ for all finite groups $H$. 
The source of the error lies in the statements about permutations
used to prove Proposition~4.1 of arXiv versions~1 and~2 of their 
paper, which are incorrect.
This claim has been withdrawn in the text of 
subsequent arXiv and the published
versions of their paper, although it still appears in the abstract
of the published version.
\end{remark}

\subsection{Update}

Since this paper was written, Nelson and Neumann~\cite{nelson-neumann08}
have shown via topological methods that $G_2(K)$ determines the knot
$K$ up to reflection. They extend their result to $G_n(K)$, $n>2$,
using the Scott-Swarup JSJ~decomposition for groups.

\subsection{Organisation}

The paper is organised as follows. In Section~\ref{Gn(L).sec} we
describe Wada and Kelly's generalised link groups, and obtain
presentations for the generalised knot groups of the square and granny
knots. We construct our target groups in Section~\ref{targets.sec},
and establish some results that will be used in the proof of
Theorem~\ref{maintheorem}, which appears in
Section~\ref{proof.sec}. Finally, we include an appendix by David
Savitt, which contains two results on \psl\ that are needed for the
proof.

\section{Generalised link groups}
\label{Gn(L).sec}

\subsection{Definition and presentations}

The groups $G_n(L)$ may be defined in several different ways.  Wada
defines them via the closed braid form of $L$, using an action of the
braid group $B_m$ on the free group $F_m$ that is compatible with the
Markov moves.  He observes that they admit a Wirtinger-type
presentation, and gives a topological description of them as the
fundamental group of a space associated with the link.  Kelly, on the
other hand, approaches them via the link quandle and the
Wirtinger-type presentation. 
We will chiefly use Wada's 
topological description, as this leads to simpler calculations
for the groups we are interested in.

To define $G_n(L)$ topologically we glue a torus to each boundary 
component of the exterior of $L$. We do this via the map 
$f: S^1\times S^1\rightarrow S^1\times S^1$ given by 
$f(z_1,z_2)=(z_1^n,z_2)$, where $z_1$ represents the meridian and
$z_2$ the longitude. The group $G_n(L)$ is defined to be the 
fundamental group of the resulting space, and is clearly an
invariant of $L$. 

\begin{figure}[t]
\leavevmode
\begin{center}
\psfrag{xi}{$x_i$}
\psfrag{xj}{$x_j$}
\psfrag{xk}{$x_k$}
\psfrag{lefthand}{$x_k = x_j^n x_i x_j^{-n}$}
\psfrag{righthand}{$x_k = x_j^{-n} x_i x_j^{n}$}
\includegraphics{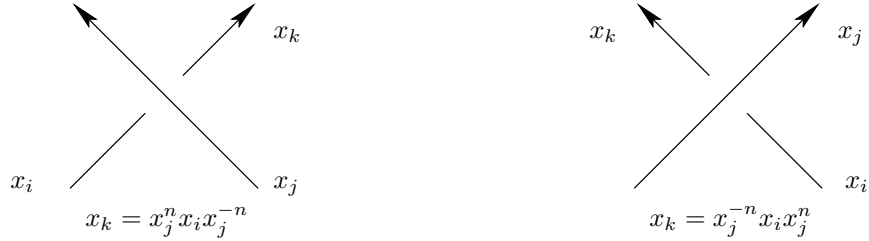}
\caption{The Wirtinger-type relations at left- and right-handed
crossings.}
\label{wirtinger.fig}
\end{center}
\end{figure}

To obtain a presentation for $G_n(L)$ we may take as generators
elements $x_i$ corresponding to each arc of the link diagram. This
leads to a Wirtinger-type presentation, with a relation 
$x_k = x_j^n x_i x_j^{-n}$ at each left-handed crossing, and
$x_k = x_j^{-n} x_i x_j^{n}$ at each right-handed crossing (see
Fig.~\ref{wirtinger.fig}). However, we may also apply the 
Seifert-van Kampen Theorem to obtain a presentation of $G_n(L)$
from a presentation 
\[
\pi_1(S^3\setminus L) = \langle g_1,\ldots, g_m | r_1,\ldots, r_p \rangle
\]
as follows. Suppose that $L$ has $\ell$ components, and let 
$\mu_i$, $\lambda_i$ be words in the generators representing a
meridian-longitude pair of the $i$th component. 
If the torus glued to the $i$th boundary component has fundamental
group $\langle \nu_i,\rho_i | \nu_i\rho_i=\rho_i\nu_i\rangle$,
then the gluing induces the identifications $\mu_i=\nu_i^n$, 
$\lambda_i = \rho_i$, and we conclude that
$G_n(L)$ has a presentation
\[
G_n(L) = \langle g_1,\ldots, g_m, \nu_1,\ldots,\nu_\ell 
| r_1,\ldots, r_p,\, \nu_i^{n}=\mu_i,\,
\lambda_i\nu_i=\nu_i\lambda_i,\,i=1,\ldots,\ell \rangle.
\]
Thus, $G_n(L)$ is obtained from the fundamental group of $L$ by adjoining
an $n$th root of each meridian that commutes 
with the corresponding longitude.
Note that, since $\mu_i$ and $\nu_i$ commute, we are not required
to use zero-framed longitudes, and may freely replace $\lambda_i$
with $\lambda_i' = \lambda_i\mu_i^k$ for any $k$.

When $K$ is a knot the presentation above reduces to
\[
G_n(K) = \langle g_1,\ldots, g_m, \nu 
| r_1,\ldots, r_p,\, \nu^{n}=\mu,\,\lambda\nu=\nu\lambda \rangle.
\]
Let $\sqrt[n]{\mu}$ denote the set of $n$th roots of $\mu$ in
$G_n(K)$, and let $C$ be the centraliser of $\mu$ in
$\pi_1(S^3\setminus K)$.  Then $C$ acts on $\sqrt[n]{\mu}$ by
conjugation, and in particular, when $K$ is composite the longitude of
each factor permutes the $n$th roots of $\mu$. This action underlies the
method by which we will  distinguish $G_n(SK)$ and $G_n(GK)$.

\begin{remark}
Crisp and Paris~\cite{crisp-paris05} generalise Wada's representations
of the braid groups and use their generalisations to define further
group invariants of links. Given a group $H$ and an element $h$ of $H$
their construction leads to an invariant $\Gamma_{(H,h)}$, such that
$G_n$ is the case $H=\integer$, $h=n$.  Topologically, their
construction replaces the attached tori with attached copies of
$X\times S^1$, where $\pi_1(X)\cong H$ and each meridian is attached
to a loop representing $h$. Careful attention is paid to the framing
in attaching the $S^1$ factor to the boundary torus.

The ability of such an invariant to distinguish links depends on the
extent to which it is able to remember the peripheral system.  Conway
and Gordon~\cite{conway-gordon75} define an extension of
$\pi_1(S^3\setminus K)$ that completely classifies knots up to oriented
equivalence. Their construction
allows them to recover the subgroups $\langle \mu\rangle$,
$\langle\lambda\rangle$ as the images of the normalisers of certain
subgroups of finite order. This method does not appear to be
directly applicable to $G_n$, since it has no torsion, but might
apply to $\Gamma_{(H,h)}$ for suitably chosen $H$ and $h$.
\end{remark}

\begin{remark}
When $n=2$ the space used to define $G_n(L)$ topologically is
a closed non-orientable $3$-manifold. Nelson and Neumann's
proof~\cite{nelson-neumann08} that $G_2(K)$ determines $K$ up
to reflection uses the JSJ~decomposition of this $3$-manifold.
\end{remark}

\subsection{The square and granny knots}

\begin{figure}[t]
\begin{center}
\leavevmode
\psfrag{a}{$a$}
\psfrag{b}{$b$}
\psfrag{c}{$c$}
\psfrag{lefthanded}{left-handed trefoil}
\psfrag{righthanded}{right-handed trefoil}
\includegraphics{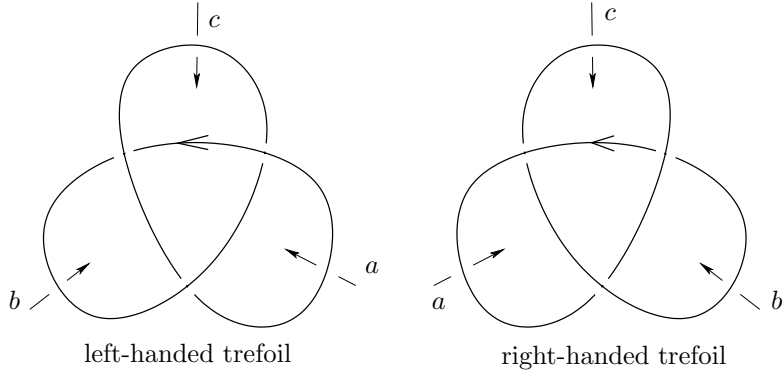}
\caption{Generators for $\pi_1(S^3\setminus K)$ for the left and
right trefoils.}
\label{trefoils.fig}
\end{center}
\end{figure}

In this section we obtain presentations for the generalised knot
groups of the square and granny knots.
Our presentations may be obtained from
those given in~\cite{lin-nelson08}, 
but our derivation uses the Seifert-van Kampen method described above,
rather than the Wirtinger-type presentation, and is more direct.

The granny knot is the connect sum of two left- or two right-handed
trefoils, while the square knot is the connect sum of a left- and
a right-handed trefoil. The fundamental groups of these knots are
therefore amalgamated free products of two copies of the trefoil 
group. Accordingly, we begin by fixing notation for the trefoil
group $T$. We choose generators for $T$ as shown in
Fig.~\ref{trefoils.fig}, so that for both the left and right
trefoils we have
\[
T \cong \langle a,b,c|ab=bc=ca\rangle
  \cong \langle a,c | aca =cac\rangle.
\]
This may be put in the form 
\[
T \cong \langle x,y| x^3=y^2\rangle
\]
by letting $x=ab$, $y=abc$, and we note that any two of
$a$, $c$, $x$ and $y$ generate $T$. The $0$-framed longitude 
corresponding to $a$ of the right-handed trefoil is represented by 
$(bac)a^{-3}$, while that of the left-handed trefoil is represented by
$(c^{-1}a^{-1}b^{-1})a^3$.
Note that $x^3= (bac)a^3$ and
$x^{-3} = (c^{-1}a^{-1}b^{-1})a^{-3}$ respectively also represent longitudes
corresponding to $a$, albeit with framings $6$ and $-6$.

\begin{figure}[t]
\begin{center}
\leavevmode
\psfrag{a}{(a)}
\psfrag{b}{(b)}
\psfrag{K1}{$K_1$}
\psfrag{m1}{$\mu_1$}
\psfrag{L1}{$\lambda_1$}
\psfrag{K2}{$K_2$}
\psfrag{m2}{$\mu_2$}
\psfrag{L2}{$\lambda_2$}
\psfrag{m=m1=m2}{$\mu=\mu_1=\mu_2$}
\psfrag{L=L1L2}{$\lambda=\lambda_1\lambda_2$}
\includegraphics{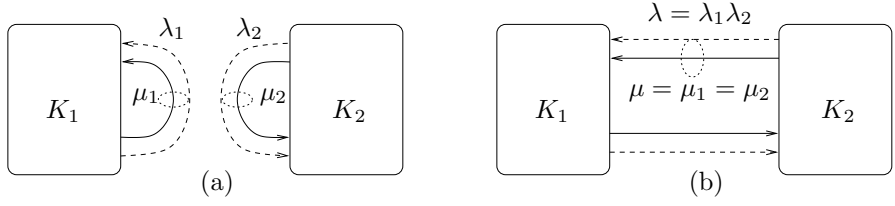}
\caption{Meridian and longitude for the connect sum of two knots.}
\label{connectsum.fig}
\end{center}
\end{figure}

To obtain presentations and longitudes of the square knot and
granny knot groups we refer to Fig.~\ref{connectsum.fig}. The group
$\pi_1(S^3\setminus(K_1\#K_2))$ is given by
\[
\pi_1(S^3\setminus(K_1\#K_2)) = 
    \pi_1(S^3\setminus K_1) *_{\langle\mu\rangle} \pi_1(S^3\setminus K_2), 
\]
and has meridian and longitude $\mu=\mu_1=\mu_2$, 
$\lambda = \lambda_1\lambda_2$.
Let 
\[
\tilde{T} = \langle d,e,f | de = ef = fd \rangle
\] 
be a second copy of $T$, with $w = de$, $z= def$, and identify the
meridians $a$ and $d$ to get
\[
\pi_1(SK)\cong \pi_1(GK) \cong \langle a,c,f | aca = cac, afa = faf\rangle.
\] 
Letting $K_1$ and $K_2$ be right- and left-handed trefoils respectively
we may take $x^3 w^{-3} = (ca)^3(fa)^{-3}$ as a longitude of $SK$;
and similarly, letting $K_1$ and $K_2$ be two right-handed trefoils, we
may take $x^3w^3=(ca)^3(fa)^{3}$ as a longitude of $GK$.
It now follows that
\begin{align*}
G_n(SK) & \cong
\langle a,c,f,\nu | aca = cac,\, afa = faf,\, \nu^n = a,\, 
         x^3w^{-3}\nu = \nu x^3w^{-3}\rangle, \\
G_n(GK) & \cong
\langle a,c,f,\nu | aca = cac,\, afa = faf,\, \nu^n = a,\, 
         x^3w^{3}\nu = \nu x^3w^{3}\rangle,
\end{align*}
which we may re-write as
\begin{align*}
G_n(SK) & \cong
\langle a,c,f,\nu | aca = cac,\, afa = faf,\, \nu^n = a,\, 
         w^{-3}\nu w^3 = x^{-3} \nu x^3\rangle, \\
G_n(GK) & \cong
\langle a,c,f,\nu | aca = cac,\, afa = faf,\, \nu^n = a,\, 
         w^{3}\nu w^{-3}= x^{-3}\nu x^3\rangle.
\end{align*}
We note that the only difference between the two presentations is
the last relation, which relates the actions by conjugation
of the longitudes of each factor on the $n$th roots of $a$.

\section{The target groups}
\label{targets.sec}

\subsection{Construction}

We will use as target groups wreath products of the form
\[
\Hpqr = \Dqr \wr \psl
      = (\Dqr)^{\pone}\rtimes \psl,
\]
where $p$, $q$, $r$ are distinct primes, \pone\ is the projective
line over the $p$-element field, and $\Dqr$ is 
a semidirect product
\[
\Dqr = (\Z{q})^{r-1}\rtimes (\Z{r}).
\]
We define these groups and some associated homomorphisms below.

To construct $\Dqr$
we regard $V = (\Z{q})^{r-1}$ as the additive group of
the finite field $\F_{q^{r-1}}$. The multiplicative group
$\F_{q^{r-1}}^\times$ is cyclic of order $q^{r-1}-1$, and so
contains an element $\zeta$ of order $r$, because
$r$ divides $q^{r-1}-1$ by Fermat's Theorem.
We may therefore define multiplication in 
$\Dqr=V\rtimes\Z{r}$ by
\[
(v,i)\cdot (w,j) = (v+\zeta^iw,i+j).
\]
We remark that when $r=2$ we have $V=\Z{q}$, and $\zeta = -1$ of
(multiplicative) order two in $(\Z{q})^\times$.  Thus $\DD{q}{2}$ is
isomorphic to $\D{q}$, the dihedral group with $2q$ elements. The
group $\DD{2}{3}$ is also isomorphic to a familiar group, the
alternating group $A_4$. Under an isomorphism $\DD{2}{3}\rightarrow
A_4$ the normal subgroup $V=(\Z{2})^2$ maps onto the Klein
$4$-group generated by the permutations of order two.

To define the wreath product we regard vectors in $(\F_p)^2$ as rows,
and use the faithful action of $\psl$ on the projective line induced
by the natural action of $\SL{p}$ on row vectors by right
multiplication.  For convenience we will identify the point in \pone\
with homogeneous co-ordinates $[x : y ]$ with the quotient
$x/y\in\Z{p}\cup\{\infty\}$; under this identification the action of
the class of $\twobytwo{a}{b}{c}{d}$ in \psl\ is given by the
fractional linear transformation
$z\mapsto\displaystyle\frac{az+c}{bz+d}$.  This right action on
indices induces a \emph{left} action of \psl\ on $(\Dqr)^{\pone}$, and
we define $\Hpqr$ to be the semidirect product of $(\Dqr)^{\pone}$ by
\psl\ defined by this action.  Elements of $\Hpqr$ then have the form
\[
\alpha=\bigl((\alpha_0,\alpha_1,\ldots,\alpha_{p-1},\alpha_\infty),\hat\alpha\bigr),
\]
where $\alpha_i \in \Dqr$ for each $i$ and $\hat\alpha\in\psl$, and 
multiplication is given by
\begin{align*}
(\alpha\beta)_i &= \alpha_i\beta_{i\cdot\hat\alpha}, &
\widehat{\alpha\beta} &= \hat\alpha\hat\beta.
\end{align*} 
Note that as permutations, $\hat\alpha$ and $\hat\beta$ are composed
from left to right.
We adopt the 
convention that $\alpha_i^{-1}$ means $(\alpha_i)^{-1}$, so that
$(\alpha^{-1})_i = \bigl(\alpha_{i\cdot\hat\alpha^{-1}}\bigr)^{-1} 
= \alpha^{-1}_{i\cdot\hat\alpha^{-1}}$.

The quotient map $\Dqr \rightarrow \Z{r}$ induces a quotient
map
\[
\Hpqr \rightarrow \Z{r}\wr \psl.
\]
We will write $[g]$ for the image of $g\in\Dqr$ in \Z{r},
so that the map $\Hpqr \rightarrow \mbox{$\Z{r}\wr \psl$}$ is given
by
\[
[\alpha] = \bigl(([\alpha_0],[\alpha_1],\ldots,[\alpha_{p-1}],[\alpha_\infty]),
                                            \hat\alpha\bigr).
\]
This map splits, and it will be convenient to distinguish a subgroup
of $\Hpqr$ isomorphic to $\Z{r}\wr\psl$. Fixing $\xi\in\Dqr$ such
that $[\xi]=1\in\Z{r}$ we let 
$\Apr=\langle\xi\rangle\wr\psl\subseteq\Hpqr$.

Since $\Z{r}$ is abelian we may quotient further to get a well defined map
\[
\Hpqr\rightarrow \Z{r}\wr \psl \rightarrow \Z{r}, 
\]
given by 
\[
[[\alpha]] = [\alpha_0]+[\alpha_1]+\cdots +[\alpha_{p-1}]+[\alpha_\infty].
\]
The subgroup $V\wr\psl$ has a similarly defined map
$\|\cdot\|:V\wr\psl\rightarrow V$. We will use these abelianisations
in Section~\ref{longitude.sec}.

To prove Theorem~\ref{maintheorem} we will choose $p$ co-prime to $n$, $q$
dividing $n$, and $r$ co-prime to $n$ and $|\psl|$.  These
divisibility requirements reflect the roles that the factors $\psl$,
$(\Z{q})^{r-1}$ and $\Z{r}$ play in the proof. The $\psl$ factor
ensures that there are nontrivial homomorphisms $G\rightarrow \Hpqr$;
the $(\Z{q})^{r-1}$ factor will give us many $n$th roots; and the
$\Z{r}$ factor will exhibit the nontrivial action by $x^3,w^3\in G$ on
the roots of the meridian.

\begin{remark}
When $p=5$ we have $\psl[5]\cong A_5$, which was the group used
by Fox~\cite{fox52} to show that the complements of the square and
granny knots can be distinguished by their peripheral subgroups.
\end{remark}

\subsection{Example}
\label{example.sec}

As an aid to understanding we illustrate the above constructions in
the case of \HEx. This group is not used in the proof of
Theorem~\ref{maintheorem} (distinguishing $w^{\pm 3}$ will ultimately
hinge on being able to distinguish $\pm 3k$ in $\Z{r}$), but has the
advantage of being a small example with a simple geometric
interpretation.

\begin{figure}[t]
\begin{center}
\psfrag{0}{$0$}
\psfrag{1}{$1$}
\psfrag{2}{$2$}
\psfrag{3}{$3$}
\psfrag{4}{$4$}
\psfrag{rho}{$\rho$}
\psfrag{sigma}{$\sigma$}
\includegraphics[scale=0.66]{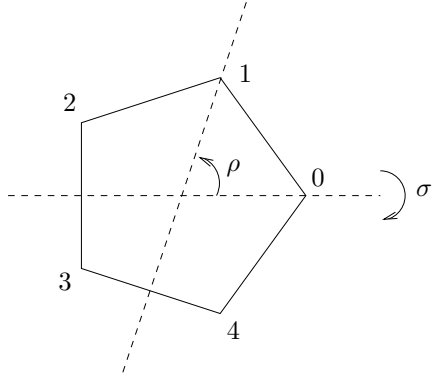}
\caption{The vertex labelling of the pentagon, and generators for $\D{5}$.}
\label{pentagon.fig}
\end{center}
\end{figure}

The permutation representation of $\psl[3]$ on $\pone[3]$ induces an
isomorphism between $\psl[3]$ and $A_4$.  Since $\DD{5}{2}$ is
isomorphic to \D{5}, we therefore have $\HEx\cong \D{5}\wr A_4$.  To
give this group a geometric interpretation let $\pent$ be the group of
symmetries of four regular pentagons.  Label the vertices of each
pentagon as in Fig.~\ref{pentagon.fig}, and let
\begin{align*}
\rho &= (0\;1\;2\;3\;4), &  \sigma &= (0)(1\;4)(2\;3)
\end{align*}
be the indicated generators of $\D{5}$. The labelling allows us to
regard maps between distinct pentagons as elements of $\D{5}$, and
induces an isomorphism $\pent\rightarrow\D{5}\wr S_4$. Our group
\HEx\ may therefore be viewed as the subgroup of $\pent$ such
that the underlying permutation of the polygons is even.

\begin{figure}[t]
\begin{center}
\psfrag{0}{$0$}
\psfrag{1}{$1$}
\psfrag{2}{$2$}
\psfrag{inf}{$\infty$}
\psfrag{alpha}{$\alpha$}
\psfrag{beta}{$\beta$}
\psfrag{gamma}{$\gamma=\alpha\beta$}
\psfrag{alphainv}{$\alpha^{-1}$}
\psfrag{a0}{$\alpha_0$}
\psfrag{a1}{$\alpha_1$}
\psfrag{a2}{$\alpha_2$}
\psfrag{ai}{$\alpha_\infty$}
\psfrag{b0}{$\beta_0$}
\psfrag{b1}{$\beta_1$}
\psfrag{b2}{$\beta_2$}
\psfrag{bi}{$\beta_\infty$}
\psfrag{g0}{$\gamma_0=\alpha_0\beta_1$}
\psfrag{g1}{$\gamma_1=\alpha_1\beta_2$}
\psfrag{g2}{$\gamma_2=\alpha_2\beta_0$}
\psfrag{gi}{$\gamma_\infty=\alpha_\infty\beta_\infty$}
\psfrag{a0inv}{$\alpha_0^{-1}$}
\psfrag{a1inv}{$\alpha_1^{-1}$}
\psfrag{a2inv}{$\alpha_2^{-1}$}
\psfrag{aiinv}{$\alpha_\infty^{-1}$}
\includegraphics[scale=0.66]{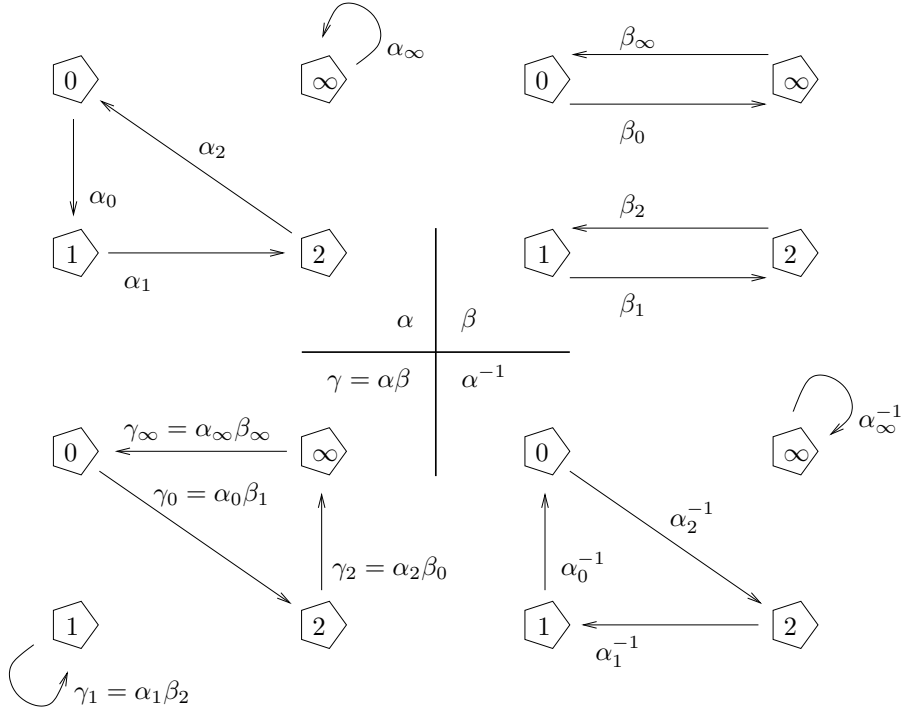}
\caption{Visualising elements and group operations in \HEx.}
\label{groupops.fig}
\end{center}
\end{figure}

To illustrate the group operations in \HEx\ define 
$A$, $B$, $C$ in $\SL{\integer}$ by
\begin{align*}
A &= \twobytwo{1}{0}{1}{1}, & B &= \twobytwo{0}{-1}{1}{0}, &
C &= AB = \twobytwo{0}{-1}{1}{-1}, 
\end{align*}
and let $\hat\alpha$, $\hat\beta$, $\hat\gamma$ be the corresponding
projections to $\psl[3]$. As fractional linear transformations we
have 
\begin{align*}
z\cdot\hat\alpha &= z+1, & z\cdot \hat\beta &= -1/z, &
z\cdot\hat\gamma &= (z\cdot\hat\alpha)\cdot\hat\beta=-1/(z+1),
\end{align*}
giving
\begin{align*}
\hat\alpha &= (0\;1\;2)(\infty), &
\hat\beta  &= (0\;\infty)(1\;2), &
\hat\gamma &= \hat\alpha\hat\beta = (0\;2\;\infty)(1)
\end{align*}
as permutations of $\{0,1,2,\infty\}$. 
Elements $\alpha,\beta\in\HEx$ permuting the pentagons according to 
$\hat\alpha,\hat\beta$
may therefore be visualised as in Fig.~\ref{groupops.fig}, where
each arrow and label describes an isometry of pentagons.  The product
$\gamma = \alpha\beta$ is found by composing $\alpha$ and $\beta$ from
left to right, and inverses are found by reversing all arrows and inverting
all labels. For example, if
\begin{align*}
\alpha &= \bigl((\rho,\rho^2\sigma,\rho^3,\rho^4\sigma),\hat\alpha\bigr), &
\beta &= \bigl((\rho^3,\rho,\rho^4,\sigma),\hat\beta\bigr),
\end{align*}
then
\begin{align*}
\gamma &= \bigl((\rho\cdot\rho,\rho^2\sigma\cdot\rho^4,\rho^3\cdot\rho^3,
\rho^4\sigma\cdot\sigma),\hat\alpha\hat\beta\bigr)
        = \bigl((\rho^2,\rho^3\sigma,\rho,\rho^4),\hat\gamma\bigr), \\
\alpha^{-1} &= \bigl((\rho^{-3},\rho^{-1},(\rho^2\sigma)^{-1},
                       (\rho^4\sigma)^{-1}),\hat\alpha^{-1}\bigr)
        = \bigl((\rho^2,\rho^4,\rho^2\sigma,\rho^4\sigma),
                        \hat\alpha^{-1}\bigr).
\end{align*}

The homomorphism $[\,\cdot\,]:\HEx\rightarrow\Z{2}\wr\psl[3]$
records whether each map $\alpha_i$ preserves or reverses orientation,
and $[[\,\cdot\,]]$ records a ``net change of orientation''. 
With
$\alpha$, $\beta$ and $\gamma$ as above we have
\begin{align*}
[\alpha] &= \bigl((0,1,0,1),\hat\alpha\bigr), &
[\beta] &= \bigl((0,0,0,1),\hat\beta\bigr), &
[\gamma] &= \bigl((0,1,0,0),\hat\gamma\bigr), \\
[[\alpha]] &= 0, & [[\beta]] &= 1, & [[\gamma]] &= [[\alpha]]+[[\beta]] = 1.
\end{align*}
The homomorphism $\|\cdot\|$ mapping
$V\wr\psl[3]\cong\langle\rho\rangle\wr A_4$ onto $\langle\rho\rangle$
is similar to 
$[[\,\cdot\,]]$ and records a ``net rotation'': if
$\delta = \bigl((\rho^3,\rho,\id,\rho^3),\hat\delta\bigr)$
then
$\|\delta\|= \rho^{3+1+0+3}=\rho^2$.
Finally, to split $[\,\cdot\,]$ we may choose say $\xi = \sigma$; then 
$\A{3}{2}=\langle\sigma\rangle\wr \psl[3]$ is a subgroup of \HEx\ mapped
isomorphically onto $\Z{2}\wr\psl[3]$ by $[\,\cdot\,]$.

Looking ahead to Section~\ref{cycleproduct.sec}, observe that the
isomorphism $\pent\rightarrow\D{5}\wr S_4$ depends on the choice of
identification of each pentagon with the labelled pentagon of
Fig.~\ref{pentagon.fig}. From the perspective of $\HEx$, changing
these identifications conjugates the group by an element $\delta$ 
such that $\hat\delta=\id$.
For example, relabelling pentagons one and two
so that $\alpha_0$ and $\alpha_1$ are the ``identity'' conjugates
$\HEx$ by $\bigl((\id,\alpha_0,\alpha_0\alpha_1,\id),\id\bigr)$, and
takes $\alpha$ to
$\bigl((\id,\id,\alpha_0\alpha_1\alpha_2,\alpha_\infty),\hat\alpha\bigr)$.
The composition $\alpha_0\alpha_1\alpha_2$ occurring here is a map from
the zeroth pentagon to itself, and as such is well defined up to
conjugacy under such relabellings. We will refer to such a composition
around a cycle of $\hat\alpha$ as a \emph{cycle-product}.

Let $\eps$ be an element of $\HEx$ such that $\hat\eps=\hat\alpha=(0\; 1\;2)$
and $\eps_{i\cdot\hat\eps}=\eps_i$ for all $i$. In other words,
$\eps$ has the form 
$\bigl((\eps_0,\eps_0,\eps_0,\eps_\infty),\hat\alpha\bigr)$.
Using a similar relabelling  to the above we
find that $\eps$ is conjugate to
$\bigl((\id,\id,\eps_0^3,\eps_\infty),\hat\alpha\bigr)$. Since the
cycle-product $\alpha_0\alpha_1\alpha_2$ necessarily has a cube
root $\eta$ in \D{5}, we may set $\eps_0=\eta$,
$\eps_\infty=\alpha_\infty$ and conclude that $\alpha$ is conjugate to
$\eps$.  
The benefit of this is that certain calculations are simpler for
$\eps$ than they are for $\alpha$: for example, $\eps^k$ is simply
$\bigl((\eps_0^k,\eps_0^k,\eps_0^k,\eps_\infty^k),\hat\alpha^k\bigr)$.

In Section~\ref{cycleproduct.sec} we will say that an element $\eps$
such that $\eps_i$ is constant on orbits of $\hat\eps$ is in
\emph{standard form}.  For $\alpha$, $\beta$ and $\gamma$ as above we
find that $\alpha$ is conjugate to
$\bigl((\sigma,\sigma,\sigma,\rho^4\sigma),\hat\alpha\bigr)$, and
$\gamma$ is conjugate to
$\bigl((\rho^4,\rho^3\sigma,\rho^4,\rho^4),\hat\gamma\bigr)$. Thus,
both $\alpha$ and $\gamma$ are conjugate to elements in standard form;
in fact, $\alpha$ may be conjugated further to
$\bigl((\sigma,\sigma,\sigma,\sigma),\hat\alpha\bigr)$, so $\alpha$ is
conjugate to an element of $\A{3}{2}$ in standard form.
However, $\beta$ cannot be conjugated to
standard form, because the cycle-product $\beta_0\beta_\infty=\rho^3\sigma$
has no square root in $\D{5}$. 
Our arguments in later sections will be
simplified by showing that we can typically avoid this situation, 
and restrict our attention to elements in standard form.

\subsection{The structure of \Dqr}

We will make use of the following facts about elements of $\Dqr$
and representations of $T$ in \Dqr:
\begin{lemma}
~

\begin{enumerate}
\item
If $g\in\Dqr$ then the order of $g$ is $1$, $q$, or $r$.
\label{orders.part}
\item
If $g,h\in\Dqr$ commute then either $g,h\in V$, or $g$ and $h$ belong
to the same cyclic subgroup of order $r$.
\label{commuting.part}
\item
If $g=(v,0)$ has order $q$ then the conjugacy class of $g$ is
\[
\{(\zeta^i v,0)| 0\leq i\leq r-1\},
\]
and if $g=(v,i)$ has order $r$ then the conjugacy class of $g$ is
$\{h:[h]=i\}$.
\label{conjugates.part}
\item
If $\{q,r\}\not=\{2,3\}$ then a homomorphism $\phi:T\rightarrow\Dqr$
factors through $\integer$.
\label{factors.part}
\end{enumerate}
\label{hqr.lem}
\end{lemma}

\begin{proof}
Since elements of $V$ have order $1$ or $q$, to prove 
part~\eqref{orders.part} it suffices to show that an element 
$(v,i)$ with $i\not\equiv 0$ has order $r$. This follows from
the fact that $\zeta$ satisfies the equation
\[
1+\zeta^i + \zeta^{2i} +\cdots +\zeta^{(r-1)i} = 0
\]
in $\F_{q^{r-1}}$. 

To prove part~\eqref{commuting.part}, suppose that $g = (v,i)$,
$h=(w,j)$ commute and that $g\not\in V$. We must show that 
$h\in\langle g\rangle$. The equation $gh = hg$ is
\[
(v+ \zeta^i w,i+j) = ( w + \zeta^j v, i+j),
\]
which gives $(1-\zeta^j)v = (1-\zeta^i)w$. Our assumption that
$g\not\in V$ is equivalent to $i\not\equiv 0 \bmod r$, which 
implies firstly that $1-\zeta^i$ is nonzero in $\F_{q^{r-1}}$, and
secondly 
that there is some $1\leq k \leq r$ such that $j\equiv ki \bmod r$.
Thus
\begin{align*}
w &= \frac{1-\zeta^{ki}}{1-\zeta^i} \,v \\
  &= \frac{(1-\zeta^{i})(1+\zeta^i+\zeta^{2i}+\cdots+\zeta^{(k-1)i})}
          {1-\zeta^i} \,v \\
  &= \bigl(1+\zeta^i+\zeta^{2i}+\cdots+\zeta^{(k-1)i}\bigr)\, v,
\end{align*}
and it follows that $h=g^k$.

Part~\eqref{commuting.part} implies that the stabiliser of $g\in\Dqr$
has order $q^{r-1}$ if $g$ has order $q$, and order $r$ if $g$ has
order $r$. Thus, the conjugacy class of $g$ has size $r$ if $g$ has
order $q$, and size $q^{r-1}$ if $g$ has order $r$.
Part~\eqref{conjugates.part} now follows from the fact that the
conjugacy class of $(v,0)$ must contain the $r$ element set consisting
of the orbit of $v$ under the action of $\langle\zeta\rangle$, and the
fact that the conjugacy class of $(v,i)$ is contained in the $q^{r-1}$
element set $\{h:[h]=i\}$.

Finally, to prove part~\eqref{factors.part} we consider the image of
$x^3 =y^2$, which generates the centre of $T$. By~\eqref{commuting.part},
if $x^3$ does not map to the identity then the image of $\phi$ lies in
an abelian subgroup, and $\phi$ must factor through the abelianisation
$T\rightarrow\integer$. On the other hand, if 
$\{q,r\}\not=\{2,3\}$ and $\phi(x)^3 = \phi(y)^2 = \id$, then at
least one of $\phi(x)$, $\phi(y)$ must be the identity, in which case
the image of $\phi$ lies in the cyclic subgroup generated by the other.
\end{proof}

\subsection{The cycle-product and some consequences}
\label{cycleproduct.sec}

In this section we gather some results on conjugacy classes,
centralisers, and $m$th powers in \Hpqr\ that will be needed in what follows.
Our main tool for proving the results will be the cycle-product, which 
may be regarded as a kind of monodromy. 

Given $\alpha\in\Hpqr$ and $i\in\pone$ we 
define the \emph{cycle-product of $\alpha$ at $i$} by 
\[
\pi_i(\alpha) = \prod_{k=0}^{\ell_i(\hat\alpha)-1} \alpha_{i\cdot\hat\alpha^k}
  =\alpha_i \alpha_{i\cdot\hat\alpha}\alpha_{i\cdot\hat\alpha^2}\cdots 
                              \alpha_{i\cdot\hat\alpha^{\ell_i(\hat\alpha)-1}},
\]
where $\ell_i(\hat\alpha)$ is the length of the disjoint cycle of
$\hat\alpha$ containing $i$. The cycle-product is thus the ordered
product, beginning at $i$, of $\alpha_j$ for $j$ in the disjoint cycle
of $\hat\alpha$ containing $i$.  We note that the value of the 
cycle-product on a given cycle depends on the starting point $i$, but
the conjugacy class does not, because
$\pi_{i\cdot\hat\alpha}(\alpha) = \alpha_i^{-1}\pi_i(\alpha)\alpha_i$.

Given $\gamma\in\Hpqr$ we will say that $\gamma$ is in 
\emph{standard form} if
\[
\gamma_{i\cdot\hat\gamma} = \gamma_i
\]
for all $i$. If in addition
$\pi_i(\gamma)=\gamma_i^{\ell_i(\hat\gamma)}=\id$ if and only if
$\gamma_i=\id$, we will say that $\gamma$ is in \emph{reduced standard
form}. Our first two results generalise the example at the end of
Section~\ref{example.sec}, and show that many elements of \Hpqr\ are
conjugate to an element in reduced standard from.

\begin{lemma}
Suppose that $\ell_i(\hat\alpha)$ is co-prime to the order of $\pi_i(\alpha)$
for all $i$. Then $\alpha$ is conjugate to an element $\gamma$ in
reduced standard form such that
$\hat\gamma=\hat\alpha$ and $\pi_i(\gamma)$ is conjugate to
$\pi_i(\alpha)$ for all $i$. 
\label{standardform.lem}
\end{lemma}

\begin{proof}
Let $\mathcal{O}$ be a set of orbit 
representatives for the action of $\langle\hat\alpha\rangle$ on \pone.
Let $\beta=\delta\alpha\delta^{-1}$, where $\delta\in\Hpqr$ is defined by 
$\hat\delta=\id$ and
\[
\delta_{j\cdot\hat\alpha^m} = \prod_{k=0}^{m-1} \alpha_{j\cdot\hat\alpha^{k}}
      = \alpha_j\alpha_{j\cdot\hat\alpha}\cdots\alpha_{j\cdot\hat\alpha^{m-1}}
\]
for each $j\in\mathcal{O}$ and $1\leq m\leq \ell_j(\hat\alpha)$.
Then $\hat\beta=\hat\alpha$, and 
\[
\beta_i = \delta_i\alpha_i\delta_{i\cdot\hat\alpha}^{-1} = 
          \begin{cases}
          \pi_i(\alpha)  & \mbox{if $i\in\mathcal{O}$}, \\
          \id            & \mbox{otherwise}.
          \end{cases}
\]

Let $j\in\mathcal{O}$. Since 
$\ell_j(\hat\alpha)$ is co-prime to the order of $\pi_j(\alpha)$, there is
$\eta_j\in\langle\pi_j(\alpha)\rangle$ such that
$\eta_j^{\ell_j(\hat\alpha)}=\pi_j(\alpha)$. Define $\gamma\in\Hpqr$
by $\hat\gamma=\hat\alpha$, and $\gamma_i = \eta_j$ if $j$ represents
the orbit of $i$. 
Then $\gamma$ is in standard form, and 
$\pi_i(\gamma)=\pi_j(\gamma)=\pi_j(\alpha)$ for all $j\in\mathcal{O}$ and
$i$ in the orbit of $j$ under $\hat\alpha$, so $\pi_i(\gamma)$ is 
conjugate to $\pi_i(\alpha)$ for all $i$.
Moreover, $\gamma$ is reduced, because 
$\gamma_i\in\langle\pi_i(\gamma)\rangle$ for all $i$.
Since $\pi_j(\gamma)=\pi_j(\alpha)$ for all $j\in\mathcal{O}$, and 
$\mathcal{O}$ is a set of orbit representatives for
$\langle\hat\gamma\rangle$ acting on \pone,
the first paragraph implies that $\gamma$ is conjugate to $\beta$, and
hence to $\alpha$.
\end{proof}

\begin{lemma}
Let $\alpha$ be an element of $\Hpqr$ in standard form such that
$\alpha_i$ has order $1$ or $r$ for each $i$. Then $\alpha$ is conjugate
to an element $\gamma$ of $\Apr$ in standard form. Moreover, if $\alpha$
is reduced then $\gamma$ may be chosen to be reduced.
\label{inHpr.lem}
\end{lemma}
Note that if $\alpha$ is in reduced standard form and $\pi_i(\alpha)$ has
order $1$ or $r$ for all $i$, then $\alpha_i$ must have order $1$ or $r$ for
all $i$ also, so $\alpha$ satisfies the hypotheses of Lemma~\ref{inHpr.lem}.
We will use the lemma in this form in Section~\ref{meridian.sec}.

\begin{proof}
By Lemma~\ref{hqr.lem},
$\alpha_i$ is conjugate to $\xi^{[\alpha_i]}$ for
each $i$. Since 
$\alpha$ is in standard form we may therefore choose $\beta_i\in\Dqr$ for
each $i$ such that
$\beta_{i\cdot\hat\alpha}=\beta_i$ and  
$\beta_i\alpha_i\beta_i^{-1}\in\langle\xi\rangle$.
Let $\gamma=\beta\alpha\beta^{-1}$, where
$\beta=\bigl((\beta_0,\beta_1,\ldots,\beta_\infty),\id\bigr)$.
Then
\[
\gamma_i = \beta_i\alpha_i\beta_{i\cdot\hat\alpha}^{-1} = 
   \beta_i\alpha_i\beta_i^{-1}=\xi^{[\alpha_i]},
\]
so $\gamma\in\Apr$. Moreover $\hat\gamma=\hat\alpha$, so 
\[
\gamma_{i\cdot\hat\gamma} 
   = \beta_{i\cdot\hat\alpha}\alpha_{i\cdot\hat\alpha}
                               \beta_{i\cdot\hat\alpha}^{-1}
   = \beta_i\alpha_i\beta_i^{-1} = \gamma_i,
\]
and $\gamma$ is still in standard form; and 
$\pi_i(\gamma)=\beta_i\pi_i(\alpha)\beta_i^{-1}$ for all $i$, so if
$\alpha$ is reduced, so is $\gamma$. 
\end{proof}

We next study certain elements of the centraliser of an element in 
reduced standard form.

\begin{lemma}
Let $\alpha$ be an element in reduced standard form, and suppose that
$\gamma$ commutes with $\alpha$. If $\alpha_i$ is constant on orbits of
$\hat\gamma$ then
$\gamma_i$ commutes with $\alpha_i$ for each $i$, and $\gamma_i$ is
constant on orbits of $\hat\alpha$.
\label{centraliser.lem}
\end{lemma}

\begin{proof}
We have
$(\alpha\gamma)_i = \alpha_i\gamma_{i\cdot\hat\alpha}$, and
$(\gamma\alpha)_i = \gamma_i\alpha_{i\cdot\hat\gamma}=\gamma_i\alpha_i$.
Thus $\alpha\gamma=\gamma\alpha$ implies 
\begin{equation}
\gamma_{i\cdot\hat\alpha}=\alpha_i^{-1}\gamma_i\alpha_i
\label{ag=ga.eq}
\end{equation}
for each $i$. 
If $\alpha_i=\id$ then we are done. Otherwise, we conjugate~\eqref{ag=ga.eq}
by $\alpha_i$ repeatedly and use the fact that 
$\alpha_{i\cdot\hat\alpha^j}=\alpha_i$ to
get
\[
\gamma_i=\gamma_{i\cdot\hat\alpha^{\ell_i(\hat\alpha)}}
   =\alpha_i^{-\ell_i(\hat\alpha)}\gamma_i\alpha_i^{\ell_i(\hat\alpha)}.
\]
Now $\alpha_i^{\ell_i(\hat\alpha)}\not=\id$, because $\alpha_i\not=\id$ and
$\alpha$ is in reduced standard form. Therefore 
$\alpha_i^{\ell_i(\hat\alpha)}$ generates $\langle\alpha_i\rangle$, because
$\alpha_i$ has prime order. It follows that $\alpha_i$ and $\gamma_i$ 
commute, and equation~\eqref{ag=ga.eq} becomes 
$\gamma_{i\cdot\hat\alpha}=\gamma_i$, as required.
\end{proof}

Finally, we give a necessary 
condition for $\alpha\in\Hpqr$ to be an $m$th power.
We will use the condition in Section~\ref{meridian.sec} with $m=n$, and
again in Section~\ref{longitude.sec}, with $m=2$ and $3$.

\begin{lemma}
\label{power.lem}
Suppose that $\alpha = \gamma^m$ in \Hpqr. Then $\hat\gamma^m=\hat\alpha$,
and
\[
\pi_i(\alpha)  = \bigl(\pi_i(\gamma)\bigr)^{m/\gcd(\ell_i(\hat\gamma),m)}.
\]
In particular, the conjugacy class of $\pi_i(\alpha)$ is constant on
orbits of $\hat\gamma$.
\end{lemma}

\begin{proof}
The equation $\alpha = \gamma^m$ gives $\hat\alpha=\hat\gamma^m$ and
\[
\alpha_i = \prod_{k=0}^{m-1} \gamma_{i\cdot\hat\gamma^k}
         = \gamma_i \gamma_{i\cdot\hat\gamma} \gamma_{i\cdot\hat\gamma^2}
                            \cdots\gamma_{i\cdot\hat\gamma^{m-1}}
\]
for all $i$. Suppose that $i$ belongs to a cycle of $\hat\gamma$ of length
$\ell$. Then $i$ belongs to a cycle of $\hat\alpha$ of length
$\ell'=\ell/\gcd(\ell,m)$, and
\[
\pi_i(\alpha) = \prod_{s=0}^{\ell'-1} \alpha_{i\cdot\hat\alpha^s} 
              = \prod_{s=0}^{\ell'-1} \alpha_{i\cdot\hat\gamma^{ms}} 
              = \prod_{s=0}^{\ell'-1} \prod_{k=0}^{m-1}
                         \gamma_{i\cdot\hat\gamma^{ms+k}}
               = \prod_{t=0}^{m\ell'-1} \gamma_{i\cdot\hat\gamma^{t}}.
\]
Let $m'=m/\gcd(\ell,m)$. Then $m\ell'-1 = \mathrm{lcm}(l,m)-1 = \ell m'-1$, so
\[
\prod_{t=0}^{m\ell'-1} \gamma_{i\cdot\hat\gamma^{t}}
              = \prod_{s=0}^{m'-1} \prod_{k=0}^{\ell-1}
                         \gamma_{i\cdot\hat\gamma^{\ell s+k}} 
              = \prod_{s=0}^{m'-1} \pi_i(\gamma) 
              = \bigl(\pi_i(\gamma)\bigr)^{m'},
\]
giving $\pi_i(\alpha) = \bigl(\pi_i(\gamma)\bigr)^{m/\gcd(\ell,m)}$, as 
required. The last statement follows from the fact that 
$\ell_i(\hat\gamma)$ and 
the conjugacy class of $\pi_i(\gamma)$ are
constant on orbits of $\hat\gamma$.
\end{proof}

\section{Proof of Theorem~\ref{maintheorem}}
\label{proof.sec}

\subsection{Strategy}
\label{strategy.sec}

As observed by Lin and Nelson~\cite{lin-nelson08} in the case
$K=SK$ or $GK$, every homomorphism
from $G_n(K)$ to a group $H$ arises as a compatible choice of a homomorphism
from $G_1(K)$ to $H$ and an $n$th root of the image of the meridian.
Given a map $\rho:G_1(K)\rightarrow H$ and an element $\eta$ of $H$,
we will say that $(\rho,\eta)$ is
a \emph{map-root pair} for $K$ in $H$ if
$\eta^n=\rho(\mu)$. Clearly, a map-root pair defines a homomorphism
$\tilde\rho:G_n(K)\rightarrow H$ precisely when it satisfies the
compatibility condition
\[
\rho(\lambda)\eta = \eta\rho(\lambda).
\] 

Since $G_1(SK)$ and $G_1(GK)$ are both isomorphic to 
\[
G = \langle a,c,f | aca = cac, \, afa = faf \rangle, 
\]
with common meridian $a$, map-root pairs for $GK$ and $SK$ in 
$\Hpqr$ co-incide. A pair $(\rho,\eta)$ is compatible for $SK$ if
\[
\rho(x^3w^{-3})\eta = \eta\rho(x^3w^{-3}),
\]
and it is compatible for $GK$ if
\[
\rho(x^3w^{3})\eta = \eta\rho(x^3w^{3}).
\]
We may regard a homomorphism
$\rho:G\rightarrow\Hpqr$ as a pair of homomorphisms
$\rho_c,\rho_f:T\rightarrow\Hpqr$ such that $\rho_c(a)=\rho_f(a)$,
using the fact that
\[
G = \langle a,c\rangle *_{\langle a\rangle} \langle a,f\rangle
  \cong T *_{\langle a\rangle} T.
\]
Given $\alpha\in\Hpqr$, let
\[
\C = \{\beta\in C_{\Hpqr}(\alpha) : \hat\beta=\id\}.
\]
We define an action of \C\ on map-root pairs 
$(\rho,\eta) = ((\rho_c,\rho_f),\eta)$ such that $\rho(a)=\alpha$ by
$\beta\cdot(\rho,\eta)=(\beta\cdot\rho,\eta)$, where
\[
\beta\cdot(\rho_c,\rho_f) = (\rho_c,\beta\rho_f\beta^{-1}).
\]
We prove Theorem~\ref{maintheorem} by showing that, for suitably
chosen $p$, $q$ and $r$,
\begin{enumerate}
\renewcommand{\theenumi}{\Roman{enumi}}
\item
for any map-root pair $(\rho,\eta)$ for $SK$ and $GK$ in \Hpqr, 
the orbit of $(\rho,\eta)$ under
the action of \C[\rho(a)]\ contains at least as many compatible
pairs for $SK$ as it does for $GK$; and
\label{inequality.part}
\item
there are compatible map-root pairs $(\rho,\eta)$ for $SK$ in \Hpqr\
such that the orbit of $(\rho,\eta)$ under the action of $\C[\rho(a)]$ 
contains no pairs that are compatible for $GK$.
\label{strict.part} 
\end{enumerate}

The prime
$r$ will be chosen so that $r$ is co-prime to $n$ and $|\psl|$. 
This implies that $r$ is co-prime to $6$, because $6$ divides the order
of \psl\ for all $p$. The primes
$p$ and $q$ will be chosen depending on $n$ as follows:
\begin{enumerate}
\item
If $n$ is not divisible by $30$ then we let $p$ be the least prime co-prime
to $n$, and let $q$ be any prime dividing $n$.
In this case $p=2$, $3$ or $5$, $\psl$ is isomorphic to 
$S_3$, $A_4$ or $A_5$, and $r$ is
co-prime to $6n$ if $p=2$ or $3$, or to $30n$ if $p=5$.
\item
If $30$ divides $n$ then we may choose $q\geq 5$ dividing $n$, and we choose
$p$ co-prime to $n$ such that \psl\ has no elements of orders $4$, $5$,
$9$, or $q$. 
\label{30|n.case}
\end{enumerate}
For $p$ co-prime to $2$ we have 
\[
|\psl| = \frac{1}{2}p(p-1)(p+1),
\]
so $\psl$ has no elements of orders $5$, $9$ or $q$ if
$p\not\equiv 0$ or $\pm 1$ mod $5$, $9$ or $q$. To 
show that we may avoid elements of order
$4$ we use the fact that $-1$ is a sum of two squares mod $p$. 
If $a^2+b^2=-1$ in \Z{p} then the matrices
\[
\twobytwo{0}{1}{-1}{0}, \quad \twobytwo{a}{b}{b}{-a}
\]
belong to \SL{p}, and for $p\geq 3$ their images in \psl\ generate a
subgroup $K$ isomorphic to
$\Z{2}\times\Z{2}$. If $p$ is congruent to
$3$ or $5$ mod $8$ then $|\psl|\equiv 4\bmod 8$, so $K$
is a Sylow 2-subgroup and there are no elements of order $4$.
Sufficient conditions on $p$ in case~\eqref{30|n.case} may thus be
expressed in terms of congruence conditions modulo $5$, $8$, $9$ and
$q$, and Dirichlet's Theorem guarantees that we may choose $p$ as
required.

The conditions on $p$, $q$ and $r$ are chosen in part so that the
order of a nontrivial $n$th power in \psl\ satisfies the following
lemma. In particular, the lemma shows that when $n$ is divisible by 
$30$ we may use Proposition~\ref{transitive.prop}. This will be
important when we study the image of the longitude in
Section~\ref{longitude.sec}.

\begin{lemma}
Let $\tau^n=\phi\not=\id$ in \psl. Then the orders of $\tau$ and
$\phi$ are co-prime to both $q$ and $r$. If $30\nmid n$ then 
$p=2$, $3$ or $5$, and
\[
\ord(\tau)=\ord(\phi) = \begin{cases}
                      \mbox{$2$ or $3$} & \mbox{if $p=2$ and $\gcd(3,n)=1$}, \\
                      p             & \mbox{otherwise}.
                      \end{cases}
\]
If $30\mid n$ then $\ord(\phi)>6$.
\label{orders.lem}
\end{lemma}

\begin{proof}
The fact that $\ord(\tau)$ and $\ord(\phi)$ are co-prime to $r$ is immediate
from the fact that $r$ is co-prime to $|\psl|$; the fact that they are co-prime
to $q$ comes from the fact that $\psl$ has no elements of order $q$ if 
$30\mid n$, and no elements of order $q^2$ when $30\nmid n$. 

The statements about $\ord(\tau)$ and $\ord(\phi)$ when $30\nmid n$
are proved by explicitly considering the orders of elements in $\psl$.
If $p=5$ then an element of $\psl[5]\cong A_5$ has order $1$, $2$, $3$ or 
$5$. Since we use $p=5$ exactly when $n$ is divisible by $6$ but not by $5$, 
a nontrivial $n$th power and its root must both have order $5$. Similar 
arguments apply when $p=2$ or $3$. 

Finally, if $30\mid n$ then $n$ is divisible by $6$, so $\phi$ is a sixth
power. In this case $2\leq \ord(\phi)\leq 6$ is impossible, because
otherwise some power of $\tau$ would have order 4, 5 or 9, contrary to
the choice of $p$.
\end{proof}

\subsection{Overview}
\label{overview.sec}

In order to prove statements~\eqref{inequality.part} 
and~\eqref{strict.part} we must calculate the $n$th roots of the
image of the meridian, and the possible images of the
longitudes $x^3$, $w^3$ of the two factor knots. We do this in
Sections~\ref{meridian.sec} and~\ref{longitude.sec} respectively,
and then prove statements~\eqref{inequality.part} and~\eqref{strict.part}
in Section~\ref{theorem.sec}. 
In Sections~\ref{meridian.sec} and~\ref{longitude.sec} we will restrict
our attention to $n$th powers in \Hpqr\ that project to nontrivial 
elements of \psl, as $n$th powers that project to the identity are
easily handled separately. The calculations
will be simplified by first showing that up to conjugacy, 
such an $n$th power is an element of \Apr\ in reduced standard form.
We assume throughout that $p$, $q$ and $r$ are chosen as described in
the preceding section.

The results and arguments of the following sections are fairly
technical, and it will be helpful to have an overview of the ideas
underlying the proof. These are most readily seen when 
$\rho(a)=\alpha=\eta^n$
is in reduced standard form, and the order of $\hat\eta$ is co-prime
to $n$. We outline the argument in this case. 
Under these conditions $\hat\eta$ is a power of $\hat\alpha$, so the
orbits of $\hat\alpha$ and $\hat\eta$ co-incide.
Lemma~\ref{centraliser.lem} thus implies that $\eta$ is in standard form, 
and we therefore have 
\[
\alpha_i = (\eta^n)_i = (\eta_i)^n
\]
in \Dqr. Since $n$ is divisible by $q$ but not by $r$, either $\alpha_i$ is
of order $r$ and $\eta_i$ is the unique $n$th root of $\alpha_i$ in 
$\la\alpha_i\ra$, or $\alpha_i=\id$ and $\eta_i\in V$. Thus, $n$th roots
$\eta$ of $\alpha$ with $\ord\hat\eta$ co-prime to $n$ are parametrised
by $V^c$, where $c$ is the number of cycles of $\hat\eta$
on which $\alpha_i=\id$. 

To calculate the possible values of $\eps=\rho(x^3)$ we will use
Theorem~\ref{thm:dichotomy}, which shows that either
$\hat\eps=\hat\alpha^6$, or $\hat\eps=\id$.  Since $\eps$ commutes
with $\alpha$, Lemma~\ref{centraliser.lem} implies that $\eps_i$ is
constant on orbits of $\hat\alpha$, and commutes with $\alpha_i$ for
all $i$. When $\hat\eps=\hat\alpha^6$ we will be able to show that in
fact $\eps=\alpha^6$, and when $\hat\eps$ is trivial
Proposition~\ref{transitive.prop} will allow us to show that
$[\eps_i]$ is constant. Consequently, in the latter case $[\eps_i]$ is 
equal to the average value $[[\eps]]/(p+1)$, which is equal to
$6[[\alpha]]/(p+1)$.  The same
arguments apply to $\delta=\rho(w^3)$, so the action on
$V^c$ induced by $\eps$ or $\delta$ acting on the roots of $\alpha$ by
conjugation is either trivial, or is given by multiplication by
$\zeta^{6[[\alpha]]/(p+1)}$.  Writing the compatibility conditions in
the form
\[
\delta^{-1}\eta\delta=\eps^{-1}\eta\eps, \qquad 
\delta\eta\delta^{-1}=\eps^{-1}\eta\eps,
\]
it is easily checked that for such roots either both, neither, or only
the left condition (that for $SK$) will be satisfied.  The strict
inequality in Theorem~\ref{maintheorem} will come from showing that this last
case can in fact occur.

In general, $\eps$ and $\delta$ will be as described above, but 
$\hat\eta$ need not be a power of $\hat\alpha$, so $\eta$
need not be in standard form. 
As a result, the set parametrising
the $n$th roots of $\alpha$ that project to a given element of \psl\ 
can be more complicated. 
The product $V^c$ will nevertheless 
still occur as a factor, with $\eps$ and $\delta$ acting as above
again playing an important role in determining compatibility. 
However, it will no longer be the 
case that a compatible pair for $GK$ is necessarily compatible for $SK$,
necessitating our use of the $\C$-action as an accounting device.

\subsection{The image and roots of the meridian}
\label{meridian.sec}

In this section we characterise up to conjugacy solutions to the pair
of equations $\eta^n=\alpha$, $\hat\alpha\not=\id$. This determines
the possible values of the image and root of the meridian in a
map-root pair $(\rho,\eta)$ with $\widehat{\rho(a)}\not=\id$.  In what
follows we will use $1/n$ to denote the multiplicative inverse of $n$
in \Z{r}.

\begin{lemma}
\label{meridian.lem}
If $\alpha$ is an $n$th power in \Hpqr\ such that $\hat\alpha\not=\id$,
then $\alpha$ is conjugate to an element of $\Apr$ in reduced standard form.
\end{lemma}
\begin{lemma}
Let $\alpha$ be an element of \Apr\ in reduced standard form such that
$\hat\alpha\not=\id$, and suppose
that $\tau$ is an $n$th root of $\hat\alpha$ in \psl. Then $\alpha$ has
an $n$th root $\eta\in\Hpqr$ such that $\hat\eta=\tau$ if and only if
$\alpha_i$ is constant on orbits of $\tau$. If so, then $\eta_i$ is constant
on orbits of $\hat\alpha$. 
Let $\sigma$ be a disjoint
cycle of $\tau$ of length $\ell$, let $d=\gcd(\ell,n)$, and let $j\in\pone$
belong to $\sigma$. Then $\eta_i$ is completely determined on $\sigma$ by
the values of $\eta_j,\eta_{i\cdot\tau},\ldots,\eta_{i\cdot\tau^{d-1}}$. The
possible values for these elements are given by the solutions to
\begin{equation}
\label{nonzero.eq}
\prod_{k=0}^{d-1} \eta_{j\cdot\tau^k} = \eta_j \eta_{j\cdot\tau} \cdots 
              \eta_{j\cdot\tau^{d-1}} = \alpha_j^{d/n}
\end{equation}
in $\langle\xi\rangle$ if $\alpha_j\not=\id$, and the solutions to 
\begin{equation}
\label{zero.eq}
\prod_{k=0}^{d-1} \eta_{j\cdot\tau^k} = \eta_j \eta_{j\cdot\tau} \cdots 
              \eta_{j\cdot\tau^{d-1}} \in V
\end{equation}
in $\Dqr$ if $\alpha_j=\id$.
\label{roots.lem}
\end{lemma}

Note that when $d=1$ (i.e., when $\tau$ is a power of $\hat\alpha$), 
Lemma~\ref{roots.lem} asserts that $\alpha$ does have $n$th roots $\eta$
such that $\hat\eta=\tau$, and that any such root is in standard form. 
Equations~\eqref{nonzero.eq} and~\eqref{zero.eq} reduce respectively to 
$\eta_j=\alpha_j^{1/n}$ and $\eta_j\in V$, in agreement with 
our treatment of this case in Section~\ref{overview.sec}.

\begin{proof}
If $\eta^n=\alpha$ then by Lemma~\ref{power.lem} we have
\[
\pi_i(\alpha)  = \bigl(\pi_i(\eta)\bigr)^{n/\gcd(n,\ell_i(\hat\eta))}.
\]
Now $q$ divides $n$ but not $\ell_i(\hat\eta)$, by
Lemma~\ref{orders.lem}, so $\pi_i(\alpha)$ is a $q$th power in
\Dqr. Therefore $\pi_i(\alpha)$ has order $1$ or $r$ for each
$i$. Since $\ell_i(\hat\alpha)$ is co-prime to $r$ for all $i$,
Lemmas~\ref{standardform.lem} and~\ref{inHpr.lem} imply that $\alpha$
is conjugate to an element of \Apr\ in reduced standard form.
This proves Lemma~\ref{meridian.lem}

We now suppose that $\alpha$ is an element of \Apr\ in reduced standard form,
and that $\eta$ is an $n$th root of $\alpha$ such that $\hat\eta=\tau$.
Let $\sigma$, $\ell$ and $d$ be as in the statement of the lemma. Then
$\sigma^n$ consists of $d$ disjoint cycles of length $\ell/d$, so if 
$i$ belongs to $\sigma$ we have
\[
\pi_i(\alpha) = \alpha_i^{\ell/d} = \xi^{\ell[\alpha_i]/d}.
\]
By Lemma~\ref{power.lem} the conjugacy class of $\pi_i(\alpha)$ is
constant on $\sigma$, so $\ell[\alpha_i]/d$ is constant on $\sigma$
by Lemma~\ref{hqr.lem}.
Since $\ell$ is co-prime to $r$ this implies that $[\alpha_i]$, and hence
$\alpha_i=\xi^{[\alpha_i]}$, is constant on $\sigma$.

We have now proved the ``only if'' direction
of Lemma~\ref{roots.lem}; we prove the ``if'' direction by solving
the equation $\eta^n=\alpha$ with $\hat\eta=\tau$
under the assumption that $\alpha_i$ is
constant on orbits of $\tau$. By Lemma~\ref{centraliser.lem} $\eta_i$
must be constant on orbits of $\hat\alpha$, and $\eta_i$ must commute
with $\alpha_i$, because $\eta$ commutes with $\eta^n=\alpha$. If $j$
belongs to $\sigma$ then $\{j,j\cdot\tau,\ldots,j\cdot\tau^{d-1}\}$ forms a 
complete set of $\hat\alpha$-orbit representatives for the orbit of
$j$ under $\sigma$, so $\eta$ is completely determined on $\sigma$ by
$\eta_j,\eta_{j\cdot\tau},\ldots,\eta_{j\cdot\tau^{d-1}}$. 

Now 
\[
\alpha_j = (\eta^n)_j = \prod_{k=0}^{n-1} \eta_{j\cdot\tau^k}
         = \left( \prod_{k=0}^{d-1} \eta_{j\cdot\tau^k} \right)^{n/d} =H^{n/d}.
\]
Suppose first that $\alpha_j\not=\id$. Then $\eta_{j\cdot\tau^k}$ must
belong to $\langle\xi\rangle$ for all $k$, by Lemma~\ref{hqr.lem},
and therefore $H$ is a solution to 
$H^{n/d}=\alpha_j$ in $\langle\xi\rangle$. This has the unique solution
$H=\alpha_j^{d/n}$, 
so $\eta_j,\eta_{j\cdot\tau},\ldots,\eta_{j\cdot\tau^{d-1}}$
are a solution to~\eqref{nonzero.eq} in $\langle\xi\rangle$, as claimed.

On the other hand, if $\alpha_j=\id$ then the condition 
$\eta_i\alpha_i=\alpha_i\eta_i$ places no restriction on $\eta_i$ on $\sigma$. 
$H$ is thus a solution to $H^{n/d}=\id$ in $\Dqr$. Now $n$ is co-prime
to $r$, by construction, and $\ell$ is co-prime to $q$, by 
Lemma~\ref{orders.lem}, so $n/d$ is divisible by $q$ but not by $r$. It follows
that $H$ must have order $1$ or $q$, so $H$ belongs to $V$.

We have now shown that if $\eta$ exists, then it must be as given in the
Lemma.
Conversely, let $\mathcal{O}$ be a set of orbit representatives for $\tau$
acting on \pone, and suppose that we are given an appropriate solution to 
equation~\eqref{nonzero.eq} or~\eqref{zero.eq} for each $j\in\mathcal{O}$.
Define $\eta_i$ on the rest of \pone\ by extending it to be constant on
orbits of $\hat\alpha$, and let 
$\eta=\bigr((\eta_0,\ldots,\eta_\infty),\tau\bigr)$. 
Then $(\eta^n)_j=\alpha_j$ for each $j\in\mathcal{O}$, and 
$(\eta^n)_{j\cdot\tau^m}$ is conjugate to $(\eta^n)_j$ by 
$\prod_{k=0}^{m-1} \eta_{j\cdot\tau^k}$ for each $j$ and $m$. The conjugating
element belongs to the centraliser of $\alpha_j$, so 
$(\eta^n)_{j\cdot\tau^m}=\alpha_j=\alpha_{j\cdot\tau^m}$, and we conclude that
$\eta^n=\alpha$. To complete the proof we note that 
equations~\eqref{nonzero.eq} and~\eqref{zero.eq} have
$r^{d-1}$ and $|\Dqr|^{d-1}\cdot|V|$ solutions respectively, so $\alpha$
does indeed have an $n$th root with shape $\tau$ under the hypothesis.
\end{proof}

\subsection{The image of the longitude}
\label{longitude.sec}

We next characterise the possibilities for the image of the longitude
$x^3$, $w^3$ of each factor knot, under the assumption that 
$\rho(a)=\alpha$ has the form considered in Section~\ref{meridian.sec}.

\begin{lemma}
\label{longitude.lem}
Let $\alpha$ be an element of $\Apr$ in reduced standard form such that
$\hat\alpha$ is nontrivial and has an $n$th root in \psl. Suppose that
$\rho:T\rightarrow\Hpqr$ is a homomorphism such that 
$\rho(a)=\alpha$, and let $\eps=\rho(x^3)$.  If $p$, $q$ and $r$ are
chosen as in Section~\ref{strategy.sec} then either 
$\eps=\alpha^6$, or 
\begin{enumerate}
\item
$\hat\eps=\id$;
\label{epshat.part}
\item
$\eps_i$ is constant on orbits of $\hat\alpha$;
\item
the conjugacy class of $\eps_i$ is constant on $\pone$;
\label{conjugacyclass.part}
\item
$\displaystyle [\eps_i] = \frac{6}{p+1}[[\alpha]]$ for all $i$; and
\label{[eps].part}
\item
$\eps_i=\xi^{6[[\alpha]]/(p+1)}$ if $\alpha_i\not=\id$.
\label{xi^[eps].part}
\end{enumerate}
\end{lemma}

\begin{proof}
We begin by noting that $\eps$ commutes with $\alpha$, because $x^3$
generates the centre of $T$.  Let $\chi=\rho(x)$, $\psi=\rho(y)$, so
that $\eps =\chi^3=\psi^2$.  Then $\hat\chi^3=\hat\psi^2=\hat\eps$ in
\psl, so by Theorem~\ref{thm:dichotomy} either $\hat\eps=\id$, or
$\hat\chi$ and $\hat\psi$ are powers of
$\hat\psi\hat\chi^{-1}=\hat\alpha$. In the latter case 
$\hat\eps=\hat\alpha^6$,
so in either case $\alpha_i$ is constant on orbits of
$\hat\eps$ and we may apply Lemma~\ref{centraliser.lem}. We conclude
that $\eps_i$ is constant on orbits of $\hat\alpha$ and commutes with
$\alpha_i$ for all $i$.  
Since $\alpha_i\in\langle\xi\rangle$ for each
$i$, it follows from Lemma~\ref{hqr.lem} that $\eps_i$ is a power of
$\xi$ whenever $\alpha_i$ is nontrivial.

The two cases now diverge, and we treat each in turn.  Suppose first
that $\hat\chi$ and $\hat\psi$ are not powers of $\hat\alpha$.  Then
$\hat\eps$ is the identity, and $\hat\chi$ and $\hat\psi$ must be of
orders $3$ and $2$ respectively, or else they would be powers of
$\hat\alpha$. Since $\hat\eps=\id$ we have $\pi_i(\eps)=\eps_i$ for
all $i$, so by Lemma~\ref{power.lem} the conjugacy class of $\eps_i$
is constant on the orbits of both $\hat\chi$ and $\hat\psi$.

We claim that in this
case $\hat\chi$ and $\hat\psi$ generate \psl,
so that the conjugacy class of $\eps_i$ is in fact constant on \pone. 
Indeed, this is
immediate from Proposition~\ref{transitive.prop} and
Lemma~\ref{orders.lem} when $30$ divides $n$.  When $30$ does not divide $n$
the group \psl\ is isomorphic to $S_3$, $A_4$ or $A_5$, and
$\langle\hat\chi,\hat\psi\rangle$ contains elements of orders $2$ and
$3$. Additional, when $p=5$ it contains an element of order $5$, by
Lemma~\ref{orders.lem}. It is easily seen that in all three cases
$\hat\chi$ and $\hat\psi$ generate $\psl$.

Since the conjugacy class of $\eps_i$ is constant on $\pone$, so is the 
value of $[\eps_i]$. To evaluate this
common value we use the abelianisation $\Hpqr\rightarrow\Z{r}$. 
On one hand we have
\[
[[\eps]] = \sum_{k\in\pone} [\eps_k] = (p+1)[\eps_i],
\]
and on the other we have 
$[[\eps]] = 6[[\alpha]]$, because $a$ generates the abelianisation of $T$.
Since $p+1$ divides $|\psl|$, which is co-prime to $r$, we may divide by 
$p+1$ mod $r$ to get $[\eps_i] = \frac{6}{p+1}[[\alpha]]$ for all $i$, as
claimed. Statement~\eqref{xi^[eps].part} now follows from this and 
the last sentence of the first paragraph.

We now turn our attention to the case where 
$\hat\chi$ and $\hat\psi$ are powers of $\hat\alpha$. 
In this case 
$\rho$ maps $T$ into the subgroup $\Dqr\wr\langle\hat\alpha\rangle$. 
If $\sigma_1,\ldots,\sigma_m$ are the disjoint cycles of $\hat\alpha$,
of lengths $\ell_1,\ldots,\ell_m$, then this subgroup may be regarded
as a subgroup of the direct product
\[
\prod_{k=1}^m \Dqr\wr\langle\sigma_k\rangle
= \prod_{k=1}^m (\Dqr)^{\ell_k}\rtimes\langle\sigma_k\rangle.
\]
Moreover, $\rho$ may be regarded as a product of maps $\rho_k$ to each
factor.  We may therefore consider each disjoint cycle separately.  In
what follows we will use $\sigma_k$ to denote both the cycle and the
set of points of \pone\ moved by this cycle.

Since
$\alpha_i$ and $\eps_i$ are constant on $\sigma_k$, the abelianisation
$[[\,\cdot\,]]: \Dqr\wr\langle\sigma_k\rangle\rightarrow\Z{r}$ gives
$[[\eps|_{\sigma_k}]]=\ell_k[\eps_i] = 6[[\alpha|_{\sigma_k}]] 
= 6\ell_k[\alpha_i]$. Hence $[\eps_i] = 6[\alpha_i]$ on $\sigma_k$, 
because $\ell_k$ is co-prime to $r$. 
If $\alpha_i$ is nontrivial on 
$\sigma_k$ then this equality implies
\[
\eps_i = \xi^{[\eps_i]} = \xi^{6[\alpha_i]} =  \alpha_i^6 = (\alpha^6)_i,
\]
where the last equality uses the fact that $\alpha$ is in standard
form.  However, if $\alpha_i=\id$ on $\sigma_k$ then we may only
conclude at this stage that $\eps_i\in V$, and we will need to work a
little harder to show that in fact $\eps_i=\id$.

By way of contradiction, suppose that $\eps_i=v\in V$
on $\sigma_k$, and that $v$ is nontrivial. We claim that $\rho_k$
maps $T$ into the subgroup $V\wr\langle\sigma_k\rangle$. 
To prove this we need only show that the image of $x$ lies in this
subgroup, because $a$ and $x$ generate $T$.
To do so we first note that $\eps$ is in reduced standard form, because 
$\ell_i(\hat\eps)$ divides $\ell_i(\hat\alpha)$ for all $i$, and 
$\ell_i(\hat\alpha)$ is co-prime to $qr$, by Lemma~\ref{orders.lem}.
Since $\eps$ and $\chi$
commute, and $\eps_i$ is constant 
on orbits of $\hat\chi$, we may
apply Lemma~\ref{centraliser.lem} to conclude that $\chi_i$ commutes
with $\eps_i$ for all $i$. On $\sigma_k$ this means that 
$\chi_i\in V$ for all $i$, by Lemma~\ref{hqr.lem}, and the claim 
is proved.

It follows that the abelianisation 
$\|\cdot\|:V\wr\langle\sigma_k\rangle\rightarrow V$ given by summing
$\beta_i$ over $\sigma_k$ is defined on 
$\rho_k(T)$. Applying this map we get
\[
\|\eps|_{\sigma_k}\| = \ell_k v = 6\|\alpha|_{\sigma_k}\| = 0.
\]
Thus $v$ is the identity afterall, because $\ell_k$ is co-prime to $q$.
Consequently $\eps_i=\id=\alpha_i^6$ on $\sigma_k$, and we are done.
\end{proof}

We refine Lemma~\ref{longitude.lem} when $[[\alpha]]=0$:

\begin{lemma}
\label{[[alpha]]=0.lem}
Let $\alpha$, $\rho$, $\eps$ and $p$, $q$, $r$ be as in
Lemma~\ref{longitude.lem}, and suppose that $[[\alpha]]=0$ but that
$\eps\not=\alpha^6$.  If $\alpha_i$ is nontrivial for some $i$ or if
$p+1$ is co-prime to $q$ then $\eps=\id$.
\end{lemma}

\begin{proof}
Since $\eps\not=\alpha^6$, $\eps$ is described by
statements~(\ref{epshat.part}--\ref{xi^[eps].part}) of
Lemma~\ref{longitude.lem}, and $\eps_i\in V$ for all $i$.   
If $\alpha_i$ is nontrivial for some $i$
then statements~\eqref{xi^[eps].part}, \eqref{conjugacyclass.part}
and~\eqref{epshat.part} imply that $\eps=\id$, and we are done. We
therefore assume that $\alpha_i=\id$ for all $i$, but that
$\eps\not=\id$.  

Since $\eps_i$ is constant on orbits of $\hat\alpha$ and
conjugate to $\eps_0$ for all $i$, there is $\beta\in\Hpqr$ such that
$\beta_i$ is constant on orbits of $\hat\alpha$,
$\beta_i\eps_i\beta_i^{-1} = \eps_0$ for all $i$, and
$\hat\beta=\id$. Define $\rho':T\rightarrow\Hpqr$ by $\rho'(g) =
\beta\rho(g)\beta^{-1}$. Then $\rho'(a) = \alpha$, but
$\rho'(x^3)=\eps'$, where $\widehat{\eps'}=\id$ and $\eps'_i=\eps_0$ for all
$i$. We now have $\eps'_i$ constant on orbits of $\widehat{\rho'(x)}$,
and the argument proceeds analogously to the
corresponding case in Lemma~\ref{longitude.lem} when
$\hat\eps=\hat\alpha^6$.  
If $\eps_0$ is nontrivial then Lemmas~\ref{centraliser.lem}
and~\ref{hqr.lem} show that $\rho'$ maps $T$ into $V\wr\psl$, and the
abelianisation gives
\[
\|\eps'\| = (p+1)\eps_0 = 6\|\alpha\| = 0
\]
in $V$. Thus $\eps_0$ is trivial afterall, because $p+1$ is nonzero
mod $q$, and we are done.
\end{proof}

Since $p\not\equiv-1\bmod q$ when $30\mid n$, the exceptions
to Lemma~\ref{[[alpha]]=0.lem} are when $\alpha_i=\id$ for all $i$ and
$p=2$, $q=3$; $p=3$, $q=2$; or $p=5$, $q=2$ or $3$. We show that in
these cases $\eps$ commutes with any $n$th root $\eta$ of $\alpha$. In
each case Lemma~\ref{orders.lem} shows that $\hat\alpha$ and $\hat\eta$
have order
$p$, and it is easily checked that $\hat\eta$ must be a power of
$\hat\alpha$. If $\alpha_i=\id$ for all $i$ then
Lemmas~\ref{roots.lem} and~\ref{longitude.lem} show that $\eta_i$ and
$\eps_i$ belong to $V$ for all $i$.  In addition, $\eps_i$ is constant
on orbits of $\hat\eta$, because these co-incide with the orbits of
$\hat\alpha$. Therefore
\[
(\eta\eps)_i = \eta_i\eps_{i\cdot\hat\eta} = \eta_i\eps_i = \eps_i\eta_i
= (\eps\eta)_i,
\]
and since $\widehat{\eta\eps} = \hat\eta = \widehat{\eps\eta}$, $\eps$ and
$\eta$ commute, as claimed.

\subsection{Proof of Theorem~\ref{maintheorem}}
\label{theorem.sec}

We now combine the results of the preceding sections and prove
Theorem~\ref{maintheorem} in the following form.

\begin{theorem}[Theorem~\ref{maintheorem}, refined]
Let $n\geq 2$ and suppose that $p$, $q$ and $r$ are chosen as in
Section~\ref{strategy.sec}.  Then
\[
|\Hom(G_n(GK),\Hpqr)| < |\Hom(G_n(SK),\Hpqr)|.
\]
\end{theorem}

\begin{proof}
Let $(\rho,\eta)$ be a map-root pair for $GK$ and $SK$ in \Hpqr. As outlined
in Section~\ref{strategy.sec}, we first show that the orbit of $(\rho,\eta)$
under the action of $\C[\rho(a)]$ contains at least as many compatible pairs
for $SK$ as it does for $GK$. In the simplest cases we will do this by
showing that a pair is compatible for $GK$ if and only if it is compatible
for $SK$, but in the most important case it will be necessary to consider
the group action.

\subsubsection{Trivial induced maps to \psl}

Suppose first that $\widehat{\rho(a)}=\id$. 
Since the conjugacy class of $a$ generates $G$, this
implies that the induced homomorphism $\hat\rho:G\rightarrow\psl$
is trivial. Thus $\rho$ may be regarded as a
homomorphism $G\rightarrow (\Dqr)^{r+1}$, and as such is a product of
homomorphisms $G\rightarrow\Dqr$. 
Since $r$ is co-prime to $6$ we have $\{q,r\}\not=\{2,3\}$, so 
by Lemma~\ref{hqr.lem} part~\eqref{factors.part} each map 
$G\rightarrow\Dqr$ factors
through $\integer$. Consequently, 
the product $\rho:G\rightarrow (\Dqr)^{r+1}$ factors
through $\integer$ also. 
We therefore have $\rho(f)=\rho(c)=\rho(a)$, so $\C[\rho(a)]$ fixes
$(\rho,\eta)$, and additionally $\rho(x^3)=\rho(w^3)=\rho(a^6)$. Since
$\eta$ commutes with $\eta^n=\rho(a)$, it commutes with both 
$\rho(x^3)$ and $\rho(w^3)$. Thus 
$(\rho,\eta)$ is a compatible pair for both $GK$ and $SK$.

\subsubsection{Nontrivial induced maps to \psl}

Now suppose that $\widehat{\rho(a)}\not=\id$. Then by
Lemma~\ref{meridian.lem} there is $\beta\in\Hpqr$ such that
$\alpha=\beta\rho(a)\beta^{-1}$ is an element of $\Apr$ in reduced
standard form. Since $\C=\beta\C[\rho(a)]\beta^{-1}$, and
$(\rho',\eta)$ is compatible for $GK$ or $SK$ if and only if
$(\beta\rho'\beta^{-1},\beta\eta\beta^{-1})$ is, it suffices to
prove~\eqref{inequality.part} under the assumption that
$\rho(a)=\alpha$. 
Let $\eps=\rho(x^3)$, $\delta=\rho(w^3)$. Then 
the compatibility conditions are
\[
\eps\delta^{-1}\eta = \eta\eps\delta^{-1} \quad\mbox{for $SK$},
\qquad
\eps\delta\eta = \eta\eps\delta \quad\mbox{for $GK$},
\]
and the possible values
for $\eps$ and $\delta$ are described by Lemma~\ref{longitude.lem}.
We consider two cases, according to whether or not at least one of
$\eps$, $\delta$ equals $\alpha^6$.

\subsubsection{Case 1: at least one of $\eps$, $\delta$ equals $\alpha^6$}

Since this property is preserved by the action of \C, we may simply show
that such a pair is compatible for $GK$ if and only if it is compatible for
$SK$. But this is immediate from the compatibility conditions, which may
be expressed in the form
\[
\delta^{\mp 1}\eta\delta^{\pm 1} = \eps^{-1}\eta\eps.
\]
At least one of $\eps$ and $\delta$ commutes with $\eta$, so each 
compatibility condition is satisfied if and only if the other commutes
with $\eta$ also.

\subsubsection{Case 2: $\eps\not=\alpha^6\not=\delta$}

If $[[\alpha]]=0$ then Lemma~\ref{[[alpha]]=0.lem} and the paragraph
that follows it show that both $\eps$ and $\delta$ commute with
$\eta$. Thus any such map-root pair is compatible for both $GK$ and
$SK$. In what follows we therefore assume that $[[\alpha]]\not=0$. By
Lemma~\ref{roots.lem} $\alpha_i$ is constant on orbits of $\hat\eta$,
and we consider separately those cycles of $\hat\eta$ on which
$\alpha_i$ is nontrivial and those where $\alpha_i=\id$.

Let $\sigma$ be a cycle of $\hat\eta$ on which $\alpha_i\not=\id$, and let
$i$ belong to $\sigma$. Then $\eta_i\in\langle\xi\rangle$, by 
Lemma~\ref{roots.lem}, and 
\[
\eps_i=\eps_{i\cdot\hat\eta}
=\delta_i=\delta_{i\cdot\hat\eta}=\xi^{6[[\alpha]]/(p+1)},
\]
by Lemma~\ref{longitude.lem}. Hence
\[
(\eta\eps\delta^{\pm 1})_i 
   = \eta_i\eps_{i\cdot\hat\eta}\delta_{i\cdot\hat\eta}^{\pm 1}
   = \eta_i\eps_{i}\delta_{i}^{\pm 1}
   = \eps_{i}\delta_{i}^{\pm 1}\eta_i = (\eps\delta^{\pm 1}\eta)_i.
\]
This shows that such cycles present no obstruction to compatibility for
either $SK$ or $GK$. 

Suppose then that $\sigma$ is a cycle of $\hat\eta$ of length $\ell$ on which 
$\alpha_i=\id$, and let $d=\gcd(\ell,n)$. 
The action of $\C$ on $(\rho,\eta)$ allows us to vary $\delta$ but not $\eps$,
so we will treat $\eps$ as fixed and regard $\delta$ as a variable to
solve for. If $i$ belongs to $\sigma$ then
compatibility for $SK$ requires 
\[
\eps_i\delta_i^{-1}\eta_i = 
               \eta_i\eps_{i\cdot\hat\eta}\delta_{i\cdot\hat\eta}^{-1},
\]
so $\delta_i$ must satisfy the recurrence relation
\begin{equation}
\delta_{i\cdot\hat\eta} = 
               \eta_i^{-1}\delta_i\eps_i^{-1}\eta_i\eps_{i\cdot\hat\eta}.
\label{SKrecurrence.eq}
\end{equation}
Given the value of $\delta_i$ as an initial condition this has the
unique solution
\begin{equation}
\delta_{i\cdot\hat\eta^m} = 
          \left( \prod_{k=0}^{m-1} \eta_{i\cdot\hat\eta^k}\right)^{-1}
          \delta_i\eps_i^{-1}
          \left( \prod_{k=0}^{m-1} 
                        \eta_{i\cdot\hat\eta^k}\right)\eps_{i\cdot\hat\eta^m}.
\label{SKsolution.eq}
\end{equation}
Similarly, compatibility for $GK$ requires that $\delta_i$ satisfies the
recurrence relation
\begin{equation}
\delta_{i\cdot\hat\eta} 
            = \eps_{i\cdot\hat\eta}^{-1}\eta_i^{-1}\eps_i\delta_i\eta_i,
\label{GKrecurrence.eq}
\end{equation}
which for each initial condition has the unique solution
\begin{equation}
\delta_{i\cdot\hat\eta^m} = \eps_{i\cdot\hat\eta^m}^{-1}
          \left( \prod_{k=0}^{m-1} \eta_{i\cdot\hat\eta^k}\right)^{-1}
          \eps_i\delta_i
          \left( \prod_{k=0}^{m-1} \eta_{i\cdot\hat\eta^k}\right).
\label{GKsolution.eq}
\end{equation}

We now consider the question of when the solutions~\eqref{SKsolution.eq}
and~\eqref{GKsolution.eq} are well defined on $\sigma$. Since $i$ and
$i\cdot\hat\eta^m$ belong to the same orbit of $\hat\alpha$ if and only if
$d\mid m$, for~\eqref{SKsolution.eq} to be well defined we must have
$\delta_{i\cdot\hat\eta^d}=\delta_i$, or 
\[
\delta_{i}\eps_i^{-1} = 
          \left( \prod_{k=0}^{d-1} \eta_{i\cdot\hat\eta^k}\right)^{-1}
          \delta_i\eps_i^{-1}
          \left( \prod_{k=0}^{d-1} \eta_{i\cdot\hat\eta^k}\right).
\]
But the product in parentheses belongs to $V$, by
Lemma~\ref{roots.lem}, and likewise $\delta_{i}\eps_i^{-1}\in V$, by
Lemma~\ref{longitude.lem}, which shows that $[\delta_{i}\eps_i^{-1}] = 0$. 
Thus any solution to~\eqref{SKrecurrence.eq} is well defined on $\sigma$. 
However, the corresponding condition for $GK$ is 
\[
\eps_i\delta_{i} = \left( \prod_{k=0}^{d-1} \eta_{i\cdot\hat\eta^k}\right)^{-1}
          \eps_i\delta_i
          \left( \prod_{k=0}^{d-1} \eta_{i\cdot\hat\eta^k}\right),
\]
and now $\eps_i\delta_{i}$ has order $r$: Lemma~\ref{longitude.lem}
gives $[\eps_i\delta_i] = 12[[\alpha]]/(p+1)$, and this is nonzero
because $r$ is co-prime to $6$. Thus, by Lemma~\ref{hqr.lem} a
solution to~\eqref{GKrecurrence.eq} is well defined on $\sigma$ exactly
when the product in parentheses is $\id$.

We next consider the action of $\C$. If $\beta\in\C$ then $\hat\beta=\id$,
and by Lemma~\ref{centraliser.lem} $\beta_i$ is constant on orbits of 
$\hat\alpha$, and belongs to $\langle\xi\rangle$ whenever $\alpha_i\not=\id$.
Conversely, it is easily checked that any such $\beta$ belongs to \C.
Thus, when $\alpha_i=\id$, the action of \C\ allows $\delta_i$ to be
chosen independently within its conjugacy class 
$\{h\in\Dqr:[h]=6[[\alpha]]/(p+1)\}$, subject only to the condition 
that $\delta_{i\cdot\hat\alpha}=\delta_i$. Consequently, every solution 
to~\eqref{SKrecurrence.eq} or~\eqref{GKrecurrence.eq} within this
conjugacy class may be realised by the action of \C\ on $\rho$.
Since any solution is completely determined by the initial condition
$\delta_i$, and~\eqref{SKsolution.eq} and~\eqref{GKsolution.eq} are
both conjugate to $\delta_i$, there are $|V|$ such solutions.

It follows from the above that compatibility of a pair $(\rho,\eta)$
for which $\eps\not=\alpha^6\not=\delta$ is completely determined by
the values of $\delta_i$ on cycles of $\hat\eta$ where $\alpha_i=\id$.
If there are no such cycles then $\beta\cdot(\rho,\eta)$ is compatible
for both $SK$ and $GK$ for all $\beta\in\C$. Otherwise, let $c$ be the
number of such cycles, and let $t$ be the number on which the
corresponding product $\prod_{k=0}^{d-1} \eta_{i\cdot\hat\eta^k}$ is
trivial. Then the above discussion shows that there are
\[
S(\rho,\eta) = [\Stab_{\C}(\delta):\Stab_{\C}(\rho)]\cdot|V|^c
\]
compatible pairs for $SK$ in the orbit of $(\rho,\eta)$, and
$S(\rho,\eta)$ compatible pairs for $GK$ if $t=c$, and none
otherwise. In each case there are at least as many compatible pairs
for $SK$ in the orbit as there are for $GK$, and we have established
statement~\eqref{inequality.part} of Section~\ref{strategy.sec}.

\subsubsection{Realisation}

We have now shown that
\[
|\Hom(G_n(GK),\Hpqr)| \leq |\Hom(G_n(SK),\Hpqr)|.
\]
To show that the inequality is strict we exhibit map-root pairs
$(\rho,\eta)$ realising the case above in which no pair in the 
orbit is compatible for $GK$. 

Define $X$, $Y$, $A$ in $\SL{\integer}$ by
\begin{align*}
X &= \twobytwo{0}{-1}{1}{1}, & Y &= \twobytwo{0}{-1}{1}{0}, &
A &= YX^{-1}=\twobytwo{1}{0}{1}{1}, 
\end{align*}
and let $\hat\chi$, $\hat\psi$, $\hat\alpha$ be the corresponding 
projections from $\SL{\integer}$ to \psl.
Then $X^3=-I=Y^2$, so $x\mapsto\hat\chi$, $y\mapsto\hat\psi$ 
define a homomorphism
$T\rightarrow\psl$ such that $a\mapsto\hat\alpha$. 
The corresponding 
fractional linear transformations are
$z\mapsto 1/(1-z)$, $z\mapsto -1/z$, and $z\mapsto z+1$, of orders
$3$, $2$ and $p$ respectively.

Define $\chi,\psi\in\Hpqr$
by $\chi_i=\xi^2$ for all $i$, and
\[
\psi_i = \begin{cases}
         \xi^4 & \mbox{if $i=0$},          \\
         \xi^2 & \mbox{if $i=\infty$},     \\
         \xi^3 & \mbox{otherwise}.         \\
         \end{cases}
\]
Then $(\chi^3)_i = \xi^6$ for all $i$, and likewise $(\psi^2)_i = \xi^6$
for all $i$, because $0$ and $\infty$ belong to the disjoint cycle
$(0\;\infty)$ of $\hat\psi$. Thus $x\mapsto\chi$, $y\mapsto\psi$
defines a homomorphism $\rho':T\rightarrow\Hpqr$, and we extend $\rho'$
to $G$ by defining $\rho'(f)=\rho'(c)$. 

Now
\[
(\rho'(a))_i = (\psi\chi^{-1})_i = \psi_i(\chi^{-1})_{i\cdot\hat\psi}
            = \psi_i\chi^{-1}_{i\cdot\hat\psi\hat\chi^{-1}}
            = \psi_i\chi_{i\cdot\hat\alpha}^{-1}=\psi_i\chi_{i+1}^{-1},
\]
so
\[
\rho'(a) = \bigl((\xi^2,\xi,\ldots,\xi,\id),\hat\alpha\bigr).
\]
The cycle-products are $\pi_i(\rho'(a))= \xi^{p+1}$ for $i\not=\infty$, and
$\pi_\infty(\rho'(a))=\id$. We may therefore conjugate $\rho'$ by a
suitably chosen $\beta$ with $\hat\beta=\id$ to get a homomorphism
$\rho:G\rightarrow\Hpqr$ such that
\[
\rho(a) = \alpha = \bigl((\xi^{(p+1)/p},\ldots,\xi^{(p+1)/p},\id),
                                            \hat\alpha\bigr)
\]
and
\[
\rho(x^3)=\rho(w^3) = \bigl((\xi^6,\ldots,\xi^6),\id\bigr).
\]
Since $n$ is co-prime to $p$ it has a multiplicative inverse $k$ mod $p$, so
by Lemma~\ref{roots.lem} 
\[
\eta_v = \bigl((\xi^{(p+1)/np},\ldots,\xi^{(p+1)/np},v),
                                            \hat\alpha^{k}\bigr)
\]
is an $n$th root of $\alpha$ for all $v\in V$. The argument above
shows that $\beta\cdot(\rho,\eta_v)$ is a compatible map-root pair for
$SK$ for all $\beta\in\C$ and $v\in V$, but is never compatible for
$GK$ unless $v=\id$. This completes the proof.
\end{proof}

\appendix

\section{Solutions to $x^3=y^2$ in \psl[\F_q], 
by David Savitt\protect\footnotemark}
\label{psl.apdx}

\footnotetext{Department of Mathematics, The University of Arizona,
         617 N.\ Santa Rita Avenue, Tucson AZ 85721, USA.
         Email: savitt@math.arizona.edu}

Let $p$ be a prime and $q$ a power of $p$.

\begin{theorem} \label{thm:dichotomy} 
Suppose $x,y \in \psl[\F_q]$ satisfy $x^3=y^2$.  Then
either $x^3=y^2=1$, or else there exists $z \in \psl[\F_q]$ such
that $x = z^2$ and $y=z^3$.
\end{theorem}

\begin{proof} First we observe that it suffices to prove the same statement
with $\Fq$ replaced throughout by an algebraic closure $\Fpbar$.
Indeed, let $G$ be the subgroup generated by $x$ and $y$. If there
exists $z \in \psl[\Fpbar]$ such that $z^2 = x$ and $z^3 = y$,
then $z = yx^{-1}$ automatically lies in $G$, hence in
$\psl[\F_q]$.  The advantage of the statement with $\Fq$ replaced
by $\Fpbar$ is that its truth is evidently unchanged if one
conjugates $G$ (and therefore $x,y$) by an element of
$\psl[\Fpbar]$.

We recall the following theorem of Dickson that classifies the
finite subgroups of $\psl[\Fpbar]$.

\begin{theorem} \cite[Secs. 255, 260]{Dickson} 
If $G$ is a finite subgroup of $\pgl{\Fpbar}$, then
one of the following holds:
\begin{enumerate}
\item[(i)] $G$ is conjugate to $\pgl{\F_{p^m}}$ or $\psl[\F_{p^m}]$ 
for some $m > 0$;
\item[(ii)] $G$ is conjugate to a subgroup of the upper triangular
matrices;
\item[(iii)] $G$ is isomorphic to $A_4$, $S_4$, $A_5$, or the dihedral
group $\D{n}$ of order $2n$ for some $n > 1$ not divisible by $p$.
\end{enumerate}
\end{theorem}

We proceed case by case.

\subsection*{Case (i)}  
By the first paragraph, we may without loss of
generality assume that $G$ is equal to $\psl[\F_{p^m}]$ or
$\pgl{\F_{p^m}}$.

Suppose first that $G = \psl[\F_{p^m}]$.  If $p^m=2,3$ we can
check the claim directly. If $p^m \ge 4$ then $\psl[\F_{p^m}]$ is
a simple group. But $H = \langle x^3 \rangle = \langle y^2 \rangle$
is clearly a normal subgroup (it is invariant under conjugation by
both $x$ and $y$, which generate $G$) and it cannot be equal to $G$.
Therefore $x^3 = y^2 =1$ in this case.

Next suppose $G = \pgl{\F_{p^m}}$.  Again if $p^m=2,3$ we check
the claim directly, while if $p^m > 3$ then the only nontrivial
normal subgroup of $\pgl{\F_{p^m}}$ is $\psl[\F_{p^m}]$ and we
can proceed as in the previous paragraph.

\subsection*{Case (ii)}  
By the first paragraph we may assume that $x,y$ are
upper triangular.  Let $X,Y$ be lifts of $x,y$ to $\gl{\Fpbar}$,
so that $X^3 \eta = Y^2$ for a nonzero constant $\eta$.
 Multiplying this equation through by $\eta^2$ gives $(\eta X)^3
= (\eta Y)^2$; replacing $X,Y$ by $\eta X, \eta Y$, we may therefore
assume that $X^3 = Y^2$.

Let $X_{11},Y_{11}$ be the upper left entries of $X,Y$ respectively.
Since $X,Y$ are upper triangular we see $X_{11}^3 = Y_{11}^2$. Since
$\Fpbar^{\times}$ is a direct limit of cyclic groups, there exists
$a \in \Fpbar^{\times}$ such that $X_{11} = a^2$ and $Y_{11} = a^3$.
We may argue similarly for the lower right entries, so that $X,Y$
have the form
$$ X = \begin{pmatrix} a^2 & \mu \\ 0 & b^2 \end{pmatrix}\,, \qquad
 Y = \begin{pmatrix} a^3 & \ \nu \\ 0 & b^3 \end{pmatrix}\,.$$
The condition that $X^3 = Y^2$ is now simply the condition that
their upper right entries are equal, which one computes to be
\begin{equation} \label{eq:upperright}
(a^2+ab+b^2)(a^2 - ab +b^2)\mu = (a+b)(a^2 - ab + b^2)\nu\,.
\end{equation} 
Observe that
$$a^6 - b^6 = (a-b)(a+b)(a^2+ab+b^2)(a^2 - ab + b^2)\,.$$ If any of
the final three factors on the right vanishes, then $a^6 = b^6$, and
moreover the equality~\eqref{eq:upperright} must become $0=0$.
Therefore $X^3, Y^2$ are both scalar and $x^3 = y^2 = 1$.

If none of the final factors on the right vanishes, then we may
define
$$ \lambda  = \frac{\mu}{a+b} = \frac{\nu}{a^2+ab+b^2} $$
and it is easy to check that
$$ X = \begin{pmatrix} a & \lambda \\ 0 & b \end{pmatrix}^2\,,\qquad
Y = \begin{pmatrix} a & \lambda \\ 0 & b \end{pmatrix}^3\,.$$ Hence
there exists $Z \in \gl{\Fpbar}$ such that $X=Z^2$ and $Y=Z^3$,
and the same is true for $x,y$.

\subsection*{Case (iii)} 
We consider each of the groups $A_4$, $S_4$, $A_5$,
$\D{n}$ in turn.

\subsubsection*{Subcase $G\cong A_4$}
The elements of $A_4$ have the following
shape:
\begin{center}
\begin{tabular}{c|c|c}
$g$ & $g^2$ & $g^3$ \\
\hline
$1$ & $1$ & $1$ \\
$(1\ 2)(3 \ 4)$ & $1$ & $(1\ 2)(3 \ 4)$ \\
$(1 \ 2 \ 3)$ & $(1 \ 3 \ 2)$ & 1
\end{tabular}
\end{center}
From the table, the only way a square may equal a cube is if both
are the identity.  Note that if we take $x$ to be any $3$-cycle and
$y$ to be any $(2,2)$-cycle then $x,y$ do indeed generate $A_4$.

\subsubsection*{Subcase $G\cong S_4$}
The elements of $S_4$ have the following
shape:
\begin{center}
\begin{tabular}{c|c|c}
$g$ & $g^2$ & $g^3$ \\
\hline
$1$ & $1$ & $1$ \\
$(1 \ 2)$ & $1$ & $(1 \ 2)$ \\
$(1\ 2)(3 \ 4)$ & $1$ & $(1\ 2)(3 \ 4)$ \\
$(1 \ 2 \ 3)$ & $(1 \ 3 \ 2)$ & 1\\
$(1 \ 2 \ 3 \ 4)$ & $(1 \ 3)(2 \ 4)$ & $(1 \ 4 \ 3 \ 2)$
\end{tabular}
\end{center}
From the table, if $x^3=y^2$ then either $x^3=y^2=1$, or else $y$ is
a $4$-cycle and $x = x^3 = y^2$ is a $(2,2)$-cycle.  In the latter
case we could take $z = y^{-1}$, but in fact the latter case is
ruled out as in this case $x,y$ would not generate $S_4$.  Note that
if $x^3=y^2=1$ then $x,y$ generate $S_4$ if and only if they are
simultaneously conjugate to $x = (1 \ 2 \ 3)$, $y = (1\ 4)$.

\subsubsection*{Subcase $G\cong A_5$} 
The elements of $A_5$ have the following shape:
\begin{center}
\begin{tabular}{c|c|c}
$g$ & $g^2$ & $g^3$ \\
\hline
$1$ & $1$ & $1$ \\
$(1\ 2)(3 \ 4)$ & $1$ & $(1\ 2)(3 \ 4)$ \\
$(1 \ 2 \ 3)$ & $(1 \ 3 \ 2)$ & 1\\
$(1 \ 2 \ 3 \ 4 \ 5)$ & $(1 \ 3 \ 5 \ 2 \ 4)$ & $(1 \ 4 \ 2 \ 5 \
3)$
\end{tabular}
\end{center}
From the table, if $x^3=y^2$ then either $x^3=y^2=1$, or else
$x,y$ are both $5$-cycles and $x = y^{-1}$.  As with $S_4$, the
latter case is ruled out by hypothesis.  Note that if $x^3 = y^2 =
1$ then $x,y$ generate $A_5$ if and only if they are
simultaneously conjugate to $y = (1 \ 2)(3 \ 4)$, $x = (1 \ 3 \
5)$.

\subsubsection*{Subcase $G\cong\D{n}$}
Recall that the group $\D{n}$ has presentation
$$\la r,s \ | \ s^2 = r^n = 1\,, sr = r^{-1} s \ra\,.$$  Since $x,y$
generate $\D{n}$ and $x^3=y^2$, it is easily checked that
$y=sr^{\ell}$ for some $\ell$ and $x^3 = y^2 = 1$.  Note that in
this case $x,y$ generate a group of order $6$, so the only
possibility for $n$ is $n=3$.   Moreover, $S_3$ is generated by
any element of order $3$ together with any element of order $2$.
\end{proof}

\begin{remark} \label{rmk:caseii} 
In case (ii) we can actually say somewhat more, namely that
as long as $a^2 - ab + b^2 \neq 0$ then we are always in the abelian
case (i.e., there exists $z$ such that $x = z^2$ and $y = z^3$).  It
remains only to check this claim when either $a+b$ or $a^2+ab+b^2$
is zero. Note that not both can be zero, or else $ab = (a+b)^2 -
(a^2+ab+b^2) = 0$, which is impossible.  If, for instance, $a+b=0$,
then $\mu=0$ in equation~\eqref{eq:upperright}.  
We may still define $\lambda
= \nu/(a^2+ab+b^2)$ and check that
$$ X = \begin{pmatrix} a & \lambda \\ 0 & b \end{pmatrix}^2\,,\qquad
Y = \begin{pmatrix} a & \lambda \\ 0 & b \end{pmatrix}^3\,.$$  The
case $a^2+ab+b^2$ is analogous.
\end{remark}

\begin{remark}  The final sentence of each subcase of case (iii)
shows that whenever one finds $S_3$, $A_4$, $S_4$, or $A_5$ as a
subgroup of $\psl[\F_q]$, it is possible to choose $x,y$ that
generate this subgroup and satisfy $x^3 = y^2 = 1$.
\end{remark}

\begin{proposition} 
\label{transitive.prop}
Suppose $x,y \in \psl[\Fp]$ satisfy $x^3 = y^2$.  If there does not
exist $z \in \psl[\Fp]$ such that $x = z^2$ and $y = z^3$, and if
$\ord(yx^{-1}) > 6$, then $\langle x,y \rangle = \psl[\Fp]$; in
particular this group acts transitively on $\mathbb{P}^1(\Fp)$.
\end{proposition}

\begin{proof} 
By Theorem~\ref{thm:dichotomy} we have $x^3 = y^2 = 1$.  Let $G =
\langle x,y \rangle$.  If we are in case (ii) of Dickson's theorem,
then by Remark~\ref{rmk:caseii} we have $a^2 - ab + b^2 = 0$.  Then
$a/b$ is a primitive $6$th root of unity (note that the hypotheses of
this Proposition imply $p\neq 2,3$) and the order of $yx^{-1}$ is $6$,
contrary to hypothesis.  Similarly, if we are in any of the subcases
of case (iii), then the possible cycle types of $x,y$ are spelled out
in the proof of Theorem~\ref{thm:dichotomy}, and one checks in each
case that the order of $yx^{-1}$ is at most $5$.

Therefore we are in case (i), and $G$ is conjugate to
$\pgl{\F_{p^m}}$ or $\psl[\F_{p^m}]$.  But $G$ is a subgroup of
$\psl[\Fp]$ and certainly cannot have order larger than the order of
$\psl[\Fp]$.  We conclude that $G = \psl[\Fp]$.
\end{proof}

\section*{Acknowledgement}

Christopher Tuffley thanks David Savitt for providing the results on \psl\ in 
the Appendix, and for other helpful discussions about \psl.


\begin{thebibliography}{1}

\bibitem{conway-gordon75}
J.~H. Conway and C.~McA. Gordon,
\newblock A group to classify knots,
\newblock {\em Bull. London Math. Soc.}, \textbf{7} (1975) 84--86.

\bibitem{crisp-paris05}
J.~Crisp and L.~Paris,
\newblock Representations of the braid group by automorphisms of groups,
  invariants of links, and {G}arside groups,
\newblock {\em Pacific J. Math.}, \textbf{221} No.\ 1 (2005) 1--27.

\bibitem{Dickson}
L.~E. Dickson, 
\newblock \emph{Linear Groups with an Exposition of
the Galois Field Theory} (Teubner, Leipzig, 1901).

\bibitem{fox52}
R.~H. Fox,
\newblock On the complementary domains of a certain pair of inequivalent knots,
\newblock {\em Nederl. Akad. Wetensch. Proc. Ser. A. {\bf 55} = Indagationes
  Math.}, \textbf{14} (1952) 37--40.

\bibitem{kelly90}
A.~J. Kelly,
\newblock Groups from link diagrams,
\newblock Ph.~D.\ Thesis, U. Warwick (1990).

\bibitem{lin-nelson08}
X.-S. {Lin} and S.~{Nelson},
\newblock {On generalized knot groups},
\newblock {\em J. Knot Theory Ramifications}, \textbf{17} No.\ 3 (2008) 
  263--272.
\newblock E-print arXiv:math.GT/0407050v4.

\bibitem{nelson-neumann08}
S.~Nelson and W.~D. Neumann,
\newblock The $2$-generalized knot group determines the knot,
\newblock arXiv:0804.0807 (2008).

\bibitem{wada92}
M.~Wada,
\newblock Group invariants of links,
\newblock {\em Topology}, \textbf{31} No.\ 2 (1992) 399--406.

\end{thebibliography}
\end{document}